\documentclass[11pt]{article}
\usepackage[letterpaper,margin=0.70in]{geometry}
\usepackage[numbers,sort&compress]{natbib}
\usepackage[T1]{fontenc}
\usepackage[utf8]{inputenc}
\usepackage{lmodern}
\usepackage{amsmath,amssymb,mathtools,amsthm}
\PassOptionsToPackage{table}{xcolor}
\usepackage{graphicx,booktabs,array,adjustbox}
\usepackage{xcolor,colortbl,siunitx}
\usepackage{tikz}
\usetikzlibrary{arrows.meta,positioning,calc,fit,backgrounds}
\usepackage[ruled,vlined]{algorithm2e}
\usepackage[section]{placeins}
\usepackage{flafter,float}
\usepackage{setspace,needspace,etoolbox}
\usepackage[font=small,labelfont=bf]{caption}
\usepackage{subcaption}
\AtBeginEnvironment{figure}{\singlespacing}
\AtBeginEnvironment{table}{\singlespacing\setlength{\tabcolsep}{5pt}\renewcommand{\arraystretch}{1.15}}
\AtBeginEnvironment{algorithm}{\small\singlespacing}
\makeatletter
\def\fps@figure{!htbp}
\def\fps@table{!htbp}
\makeatother
\pretocmd{\section}{\FloatBarrier\Needspace{6\baselineskip}}{}{}
\pretocmd{\subsection}{\Needspace{5\baselineskip}}{}{}
\usepackage[hidelinks]{hyperref}
\newcommand{\eps}{\varepsilon}
\newcommand{\dd}{\,\mathrm{d}}

\newtheorem{proposition}{Proposition}
\title{Inferring constitutive forces from noisy Allen--Cahn movies with moment equations}
\author{Soobeen Jung \qquad Hyunju Kim\thanks{Corresponding author: hjkim@kentech.ac.kr}\\[0.35em]
\small Department of Energy Engineering, Korea Institute of Energy Technology,\\
\small Naju 58330, Republic of Korea\\
\small sbeenjung@kentech.ac.kr}
\date{}
\begin{document}
\maketitle
\begin{abstract}
Inferring an Allen--Cahn constitutive force from images requires more than an
accurate optimizer: nonlinear observation noise alters the mean of force
features, shared noise couples instruments to the response, and a temporal
balance borrowed from another solver can remain biased even on clean data.
We develop observation-aware, noise-adjusted moment equations that address
these effects before inversion. The surface-tension-calibrated constitutive
transformer (SCCT) couples the resulting moments to a positive Bernstein
inverse and can adapt its prior and precision from images and moments. A
conditional stability bound separates moment residuals from prior error.
Synthetic tests show that correcting the estimating equation matters more
than increasing regularizer flexibility; the benefit of learned adaptation
is smaller and depends on the training comparison. Forces inferred from
planar movies predict unseen planar and three-dimensional evolutions. A
ring-cage breakup shows why diffuse-interface accuracy distinguishes forces
even when their topological predictions agree, while external surface
tension sets the energy scale separately.
\end{abstract}

\section{Introduction}
A phase-field movie shows how an interface moves, but not directly which
constitutive force drives it. Recovering that force would turn the observed
pattern into a model for evolution from a new initial field. In Allen--Cahn
dynamics, temporal change combines diffusion and a nonlinear force
\cite{AllenCahn1979,Cahn1958}, while measurement noise affects the two terms
differently. The inverse therefore begins by constructing an equation that
separates constitutive behavior from observation error.

Trajectory inversion has shown how full pattern images can determine constitutive
relations and their uncertainty \cite{Zhao2020Pattern,Zhao2021ImageInversion,Cohen2026Differentiable}.
Physics-informed, pseudo-spectral, weak-form, and integrable neural models provide related
routes to state or law estimation \cite{Raissi2019PINN,JiaZhao2021SPINN,Bao2025PFWNN,Teichert2019IDNN}.
Sparse identification of nonlinear dynamics (SINDy) and its weak-form extension
(WSINDy) instead aggregate the governing equation before regression
\cite{Brunton2016SINDy,Rudy2017PDE,Schaeffer2017Integral,Gurevich2019Robust,Reinbold2020Noisy,%
MessengerBortz2021Galerkin,Messenger2021WSINDy,Wang2019VSI,Stephany2024WeakPDELearn}.
Weak formulations reduce differentiation noise, and subsequent work analyzes consistency
and residual covariance \cite{MessengerBortz2024Consistency,Bortz2023WENDy,Rummel2025WENDyMLE}.
Together these approaches motivate a more specific question: which errors
must be removed from the estimating equation before an adaptive inverse can
usefully act on it?

Two measurement biases survive averaging. Substitution of a noisy phase value into a nonlinear force
changes its expectation, while instruments constructed from the same observation sequence are correlated
with response noise. Gaussian Hermite polynomials remove the first bias, following the
corrected-score principle \cite{Stefanski1989Unbiased,Nakamura1990Corrected,Carroll2006Measurement}
and the Gaussian moment inverse in WSINDy analysis
\cite[Lemma B.2]{MessengerBortz2024Consistency}. Gaussian recorruption supplies opposite
views for the second correction \cite{Pang2021R2R,Monroy2025R2R}. Together they restore an
expectation identity for unnormalized moments under known iid Gaussian noise. A third
issue is numerical: a balance based on the wrong time generator can remain biased even
on clean data. When generator metadata are known, we derive a discrete balance from those
time steps; otherwise we use an observation-time weak balance, following the integrator-consistency
principle of \cite{GoyalBenner2022RK4SINDy}.

This leads to the paper's central claim: for the tested inverse problem, a
sound estimating equation contributes more than a more elaborate regularizer.
The surface-tension-calibrated constitutive
transformer (SCCT) couples Gaussian polynomial correction, independent-view
instruments, and a balance chosen for the observation source to the same
constrained Bernstein inverse. It optionally uses
images and moments to set a trajectory-specific prior and precision without
trajectory-wise neural optimization
\cite{Afkham2021Regularization,AlbertiDeVito2021OptimalTikhonov,%
Hertrich2026LearnedRegComparison}. A conditional bound describes how the
moment residual and prior displacement pass through the realized inverse.
Retrained experiments then test these terms separately. Moment construction
produces the larger improvement; the original adaptive models show a smaller
advantage on original and fresh laws, but common-training precision and
structural controls do not reproduce it. The positive family
and external surface-tension calibration follow the companion study
\cite{JungKim2026STAC}; calibration follows force inference and changes neither
the force nor its forecasts. This separation is consistent with gauge questions
in related inversions
\cite{KianLiimatainenLin2024,BrunkEggerHabrich2023,NiLai2024Uniqueness}.

The evidence follows the inverse from measurement to prediction. A factorial
comparison isolates the two noise corrections after retraining neural and
non-neural inverses. Crossing forward generators with temporal balances then
tests whether the equations describe the observations that produced them.
Independent laws distinguish learned regularization from moment construction,
and information-loss and noise-model studies locate the limits of recovery.
Finally, forecasts from unseen initial fields test transfer; the three-dimensional
case separates field accuracy from a topology change that can persist under
inaccurate forces.

Section~\ref{sec:model} establishes the identifiable force, and
Section~\ref{sec:method} develops the estimating equation and inverse.
Section~\ref{sec:design} follows the numerical evidence through transfer,
with proofs and secondary controls in the appendices.

\section{Identifiable model and admissible family}\label{sec:model}
Before constructing an inverse, we distinguish the force determined by
trajectories from the energy scale that requires an external measurement.
The admissible family then makes that force inference finite dimensional.
\subsection{The identifiable constitutive target}
On a periodic domain $\Omega\subset\mathbb R^2$,
\begin{equation}
 u_t=q(\Delta u-G(u)),\qquad G(u)=F'(u)/\eps^2,
 \label{eq:model}
\end{equation}
with known $q>0$; all computations use $q=1$, corresponding to the time coordinate $qt$.
With mobility $M=q/\eps$, this is the $L^2$ gradient
flow of
\[
 E_\eps[u]=\int_\Omega\left(\frac{\eps}{2}|\nabla u|^2+\frac{F(u)}{\eps}\right)\dd x,
 \qquad \frac{\mathrm dE_\eps}{\mathrm dt}=-\frac{\eps}{q}\int_\Omega u_t^2\dd x.
\]
The order parameter is nonconserved, and this identity describes continuum dissipation;
energy behavior of the time-discrete generator still depends on its time step.

The energy identity also reveals what trajectories cannot determine: the transformation
$(\eps,F,M)\mapsto(a\eps,a^2F,M/a)$ leaves the dynamics unchanged for any $a>0$.
The trajectory-determined target is consequently $G$, although its accurate recovery still
depends on the observations and inverse assumptions. We infer $G$ first; an external surface
tension then sets the remaining energy scale (Section~\ref{sec:new_forecasts}).

\subsection{Admissible force and observations}
Within this identifiable target, we represent the force for $|u|\leq1$ by
five coefficients:
\begin{equation}
 G_c(u)=a_0\sum_{j=0}^{4}c_jb_j(u),\quad
 b_j(u)=-u(1-u^2)B_j^4(u^2),\quad c_j\geq c_{\min},
 \label{eq:force}
\end{equation}
where $B_j^4(z)=\binom4jz^j(1-z)^{4-j}$. The chosen scale $a_0$ keeps typical coefficients
near unity; the small positive floor $c_{\min}$ excludes degenerate wells. Their implementation
values are given in \ref{app:moments}. Positivity ensures the force's sign on either side of zero and an even double
well. This family and its admissibility properties are those of the companion study
\cite{JungKim2026STAC}; structural priors of this kind are also used to restrict weak-form
dictionaries before regression \cite{Tang2026PriorIDENT}. Its primitive is
\[
 H_c=a_0\sum_jc_jV_j,\qquad
 V_j(u)=\frac12\binom4j\int_0^{1-u^2}r^{5-j}(1-r)^j\dd r.
\]
For numerical trajectories that pass beyond the phase interval, the energy has a quadratic
extension and the force a linear extension.
The force and its first derivative are continuous at both endpoints
(\ref{app:moments}). The five positive coefficients vary the force shape
while retaining two pure-phase minima and one central unstable state. Inference is
parameter estimation within this admissible family. Asymmetric wells and displaced minima fall
outside it, while sparse sampling of intermediate phase values can leave even an in-family force
poorly determined.

To estimate the force within this family, we use observed phase values rather than only
interface contours. Specifically, $y_{n,i}=u(t_n,x_i)+\eta\zeta_{n,i}$ with independent and
identically distributed (iid) standard Gaussian $\zeta$ unless stated otherwise. The noise
level is supplied or estimated from the observations; experimental resolutions are specified in
Section~\ref{sec:design}. It is distinct from the downstream surface-tension datum $\sigma$.

\section{Noise-adjusted moment equation and inverse}\label{sec:method}
\begin{figure}[!htbp]\centering
\definecolor{flowNavy}{HTML}{17324D}
\definecolor{flowBlue}{HTML}{2F6FA3}
\definecolor{flowTeal}{HTML}{2A9D8F}
\definecolor{flowGold}{HTML}{E9A23B}
\definecolor{flowMist}{HTML}{F3F7FA}
\resizebox{0.78\textwidth}{!}{%
\begin{tikzpicture}[
  font=\scriptsize,
  flow/.style={-{Latex[length=2mm,width=1.2mm]},draw=flowNavy!85,line width=0.75pt},
  wire/.style={draw=flowNavy!65,line width=0.7pt},
  context/.style={flow,dashed,draw=flowBlue},
  stage/.style={draw=flowNavy!72,line width=0.7pt,rounded corners=2mm,
    fill=white,align=center,inner sep=4pt,minimum height=10mm},
  correction/.style={stage,draw=flowTeal!85!black,fill=flowTeal!7},
  balance/.style={stage,draw=flowBlue!85!black,fill=flowBlue!7},
  adaptive/.style={stage,draw=flowBlue!90!black,fill=flowBlue!12},
  tag/.style={stage,draw=flowNavy!45,minimum height=6mm,inner ysep=2.5pt},
  badge/.style={circle,fill=flowNavy,text=white,font=\bfseries\scriptsize,inner sep=1.6pt},
  steplabel/.style={font=\bfseries\scriptsize,text=flowNavy,anchor=west},
  panel/.style={draw=flowNavy!22,fill=flowMist,rounded corners=3mm,inner sep=6pt}
]
\node[badge] at (0.35,7.0) {1};
\node[steplabel] at (0.67,7.0) {Construct the moment equation};
\node[stage,text width=2.15cm] (obs) at (1.5,5.95)
  {\textbf{Noisy phase observations}\\$y(t_n,x)$\\{\scriptsize times, $q$, noise level}};
\node[correction,text width=2.55cm] (features) at (4.7,5.95)
  {\textbf{Noise-adjusted features}\\{\scriptsize Hermite terms; paired views}};
\node[balance,text width=2.5cm] (source) at (7.9,5.95)
  {\textbf{Choose the balance}\\{\scriptsize discrete or weak in observed time}};
\node[stage,text width=2.2cm,fill=flowGold!10,draw=flowGold!85!black]
  (moments) at (11.05,5.95)
  {\textbf{Moment system}\\$\overline{\mathsf H}\in\mathbb R^{9\times5},\ \bar r$};
\draw[flow] (obs.east)--(features.west);
\draw[flow] (features.east)--(source.west);
\draw[flow] (source.east)--(moments.west);

\node[badge] at (0.35,4.08) {2};
\node[steplabel] at (0.67,4.08) {Regularization};
\node[adaptive,text width=2.15cm] (scct) at (3.5,2.7)
  {\textbf{SCCT}\\{\scriptsize image + moments\\prior and precision}};
\node[stage,text width=2.15cm] (global) at (7.0,2.7)
  {\textbf{Global}\\{\scriptsize common prior and precision}};
\node[stage,text width=2.15cm] (smooth) at (10.5,2.7)
  {\textbf{Smoothness}\\{\scriptsize penalty chosen by validation}};
\draw[wire] (moments.south)--(11.05,4.38);
\draw[wire] (11.05,4.38)--(3.5,4.38);
\draw[flow] (3.5,4.38)--(scct.north);
\draw[flow] (7.0,4.38)--(global.north);
\draw[flow] (10.5,4.38)--(smooth.north);
\draw[context] (obs.west)--(-0.42,5.95)--(-0.42,2.7)--(scct.west);
\node[fill=white,inner sep=1pt,text=flowBlue] at (0.85,2.7) {image context};

\node[badge] at (0.35,0.35) {3};
\node[steplabel] at (0.67,0.35) {Infer and reuse};
\draw[wire] (3.5,1.33)--(10.5,1.33);
\draw[wire] (scct.south)--(3.5,1.33);
\draw[wire] (global.south)--(7.0,1.33);
\draw[wire] (smooth.south)--(10.5,1.33);
\node[fill=white,inner sep=1pt,text=flowNavy] at (5.25,1.33) {select one};
\node[stage,text width=3.0cm,fill=flowNavy!6] (inverse) at (7.0,0.27)
  {\textbf{Constrained inverse}\\$c_j\ge c_{\min}$};
\draw[flow] (7.0,1.33)--(inverse.north);
\node[stage,text width=2.05cm,fill=flowTeal!9,draw=flowTeal!85!black]
  (law) at (11.15,0.27) {\textbf{Inferred force}\\$G_{\widehat c}$};
\draw[flow] (inverse.east)--(law.west);
\node[tag,text width=1.8cm] (forecast) at (9.95,-1.5) {new-field forecast};
\node[tag,text width=2.1cm] (energy) at (12.5,-1.5)
  {external $\sigma$\\$\widehat F,\widehat\eps$};
\draw[wire] (law.south)--(11.15,-0.74)--(9.95,-0.74);
\draw[wire] (11.15,-0.74)--(12.5,-0.74);
\draw[flow] (9.95,-0.74)--(forecast.north);
\draw[flow] (12.5,-0.74)--(energy.north);
\begin{scope}[on background layer]
  \node[panel,fit=(obs)(features)(source)(moments)] {};
  \node[panel,fit=(scct)(global)(smooth)] {};
  \node[panel,fit=(inverse)(law)(forecast)(energy)] {};
\end{scope}
\end{tikzpicture}}
\caption{
  Method hierarchy. Hermite polynomials and paired views address distinct noise effects.
  The observation source determines the discrete or weak balance. Each regularizer uses
  the same moment system and positive constrained inverse; the dashed path supplies
  image context only to SCCT. External surface tension sets the energy scale after
  force inference and does not alter forecasts.
}
\label{fig:workflow}
\end{figure}
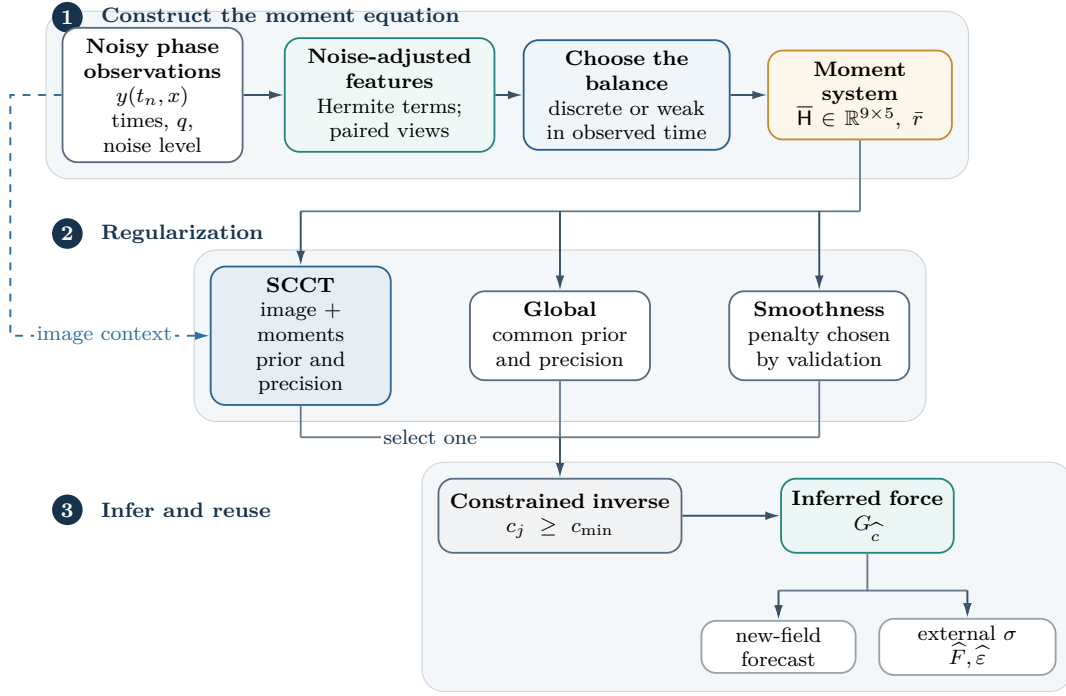

The method first removes systematic errors from the estimating equation and
then solves a common positive inverse. Figure~\ref{fig:workflow} distinguishes
that shared equation from the three regularization choices and the subsequent
energy calibration.
\subsection{Noise-adjusted estimating moments}
We begin with the nonlinear substitution bias. For Gaussian variance $v$,
initialize $Q_0,Q_1$ and, for $k\geq2$, use
\begin{equation}
 Q_0=1,\quad Q_1=y,\quad
 Q_k(y;v)=yQ_{k-1}(y;v)-(k-1)vQ_{k-2}(y;v).
 \label{eq:hermite}
\end{equation}
Since $\mathbb E Q_k(u+\xi;v)=u^k$ for $\xi\sim N(0,v)$, replacing each monomial in $b_j$ gives a
Hermite-adjusted polynomial $\mathcal B_j$. This identity is exact only when the Gaussian variance used in
$Q_k$ is the true variance. An interpolated observation uses
$v=\eta^2\sum_\ell w_\ell^2$; with a noise level estimated from the same observations, the correction is
plug-in and its unbiasedness is only approximate. Polynomial arguments are not clipped;
interpolation/model discrepancies remain.
The basis $b_j$ has degree eleven in $u$, so terms through $Q_{11}$ are needed to apply this
identity to the entire force family. This removes nonlinear-substitution bias in each feature,
but does not separate noise in the instrument from noise in the response.

For that second source of bias, form $y^\pm=y\pm\xi$ with independent
$\xi\sim N(0,\eta^2I)$ \cite{Pang2021R2R}. Under the known iid Gaussian model,
the two view errors are conditionally independent, each with variance $2\eta^2I$:
$\operatorname{Cov}(y^+,y^-\mid u)=\eta^2I-\eta^2I=0$.
They reuse one observation sequence and double its variance. Plug-in noise estimation makes the
cancellation approximate; alternative noise laws or instruments from separate trajectories need
different constructions \cite{Monroy2025R2R,KuangLin2026IV}.

The corrected features and paired views enter one estimating system. For either balance, let
$X_j$ denote a noise-adjusted basis feature and $d$ the
diffusion--time response divided by $a_0$. Temporal filtering $T$ and periodic Gaussian
smoothing $S_w$ reduce high-frequency variability. The Fourier Laplacian $L_h$ is still
applied to the data: only the observation-time weak path transfers the temporal derivative
to a test function. Nine overlapping phase instruments use
\begin{equation}
\psi_k(z)=-z(1-z^2)\exp\left[-\frac{(z^2-\mu_k)^2}{2h_\psi^2}\right],\quad k=0,\ldots,8.
\label{eq:instruments_main}
\end{equation}
Here $\mu_k$ are squared-phase centers and $h_\psi$ is their common width.
Nine moment rows overdetermine five force coefficients, although this count alone does not
guarantee full rank.
Their smoothed, clipped instrument values $Z_k$ cover $u^2\in[0,1]$; the Hermite-adjusted force
arguments remain unclipped. The two balances are described next, and
\ref{app:moments} specifies their filters, stage reconstruction, and supports.
With these ingredients, opposite views are paired using two independently
drawn artificial perturbation fields $b=1,2$:
\begin{equation}
 \begin{aligned}
 \mathsf H_{kj}&=\tfrac14\sum_{b=1}^2[\langle Z_k(y_b^+),X_j(y_b^-)\rangle+\langle Z_k(y_b^-),X_j(y_b^+)\rangle],\\
 r_k&=\tfrac14\sum_{b=1}^2[\langle Z_k(y_b^+),d(y_b^-)\rangle+\langle Z_k(y_b^-),d(y_b^+)\rangle].
\end{aligned}\label{eq:moments}
\end{equation}
Both artificial pairs share the original measurement noise, so they are not independent
replicates; their average reduces artificial-perturbation variability only.
Brackets average valid times, spatial centers, and two smoothing widths. Within a view, polynomial
correction uses variance $2\eta^2\sum_\ell w_\ell^2$. Write $X_0,d_0$ for the corresponding clean features and response. Opposite-view independence gives
\[
 \mathbb E[\langle Z^+,X^-c_\star-d^-\rangle\mid u]
 =\langle\mathbb E[Z^+\mid u],X_0c_\star-d_0\rangle.
\]
Thus, under the known iid Gaussian variance model, polynomial correction removes nonlinear
substitution bias and view pairing removes the instrument--response cross-noise term. Plug-in noise
estimation makes both cancellations approximate, and clean interpolation and model residuals remain.
The identity applies to unnormalized moments, not to the final estimator.
\ref{app:moments} traces the effects of random normalization and inversion;
the realized inverse is analyzed by the conditional bound below.

To compare rows drawn from different phase groups, normalize each moment row and its response
by the row norm, with a numerical
denominator floor specified in \ref{app:moments}. With the resulting
$\overline{\mathsf H},\bar r$, set
\[
 K=\overline{\mathsf H}^{\top}\overline{\mathsf H}/9,\qquad
 g=\overline{\mathsf H}^{\top}\bar r/9,\qquad s_K=\operatorname{tr}K/5.
\]
When the row floors are inactive, $\operatorname{tr}K=1$ and $s_K=1/5$.
The normalized spectrum measures directional conditioning rather than raw signal strength.

\subsection{Discrete and observation-time balances}\label{sec:continuum_moments}
A corrected force feature still needs a compatible temporal balance. For the reference
implicit--explicit (IMEX) generator \cite{AscherRuuthWetton1995}, diffusion acts at the later internal stages and the
force at the preceding stages. Summing those steps and reconstructing their states from
observed frames gives Eq.~\eqref{eq:imex_features_main} in \ref{app:moments}.
The stage placement preserves the generator identity; interpolation contributes a remaining
residual. This is the same consistency principle used in
Runge--Kutta-inspired regression \cite{GoyalBenner2022RK4SINDy}.

Without internal-stage metadata, use a windowed bump $\phi(t)=(1-s(t)^2)^4$, $s\in[-1,1]$ the affine time coordinate on the window.
Integration by parts gives
\begin{equation}
 \int\phi G(u)\dd t=\int\phi\Delta u\dd t+q^{-1}\int\phi'u\dd t.
 \label{eq:weak_time}
\end{equation}
For observed times $t_i$ in a window, let $\omega_i$ be trapezoidal weights, $\phi_i=\phi(t_i)$,
$\phi'_i=\phi'(t_i)$, and $D_\phi=\sum_i\omega_i\phi_i$. The actual basis and response are
\begin{equation}
 \begin{aligned}
 X_j^{\rm w}&=S_w\!\left[\frac{\sum_i\omega_i\phi_i\mathcal B_j(y_i;2\eta^2)}{D_\phi}\right],\\
 d^{\rm w}&=\frac1{a_0}S_w\!\left[\frac{L_h\sum_i\omega_i\phi_i y_i+q^{-1}\sum_i\omega_i\phi'_i
 y_i}{D_\phi}\right].
 \end{aligned}\label{eq:weak_observed}
\end{equation}
Here $y_i$ denotes one recorrupted view and $\mathcal B_j(\cdot;v)$ its Hermite-adjusted force
basis, distinct from the Bernstein polynomial $B_j^4$. Windows advance by one observed frame and use
only their observed samples; no internal states are reconstructed. Hermite polynomials use
single-frame variance $2\eta^2$ within each view, and opposite-view instruments are averaged over
the same window. The resulting system retains the form of Eq.~\eqref{eq:moments}. Validation selects
the window support to balance smoothing against temporal-quadrature error
(\ref{app:moments}). These weak moments
\cite{Messenger2021WSINDy,Wang2019VSI} require observation times, $q$, and noise level, but retain
temporal-quadrature and spatial-discretization errors.

\subsection{Positive inverse and optional learned regularization}
With noise and time-discretization effects separated, the inverse balances the
information in the moments against a
data-dependent prior. Positive inverses of Bernstein degrees $m=1,2,4$ are elevated to five
coefficients and blended as
$c_{\rm pr}=\sum_m\pi_m c^{(m)}$, with $\pi_m\geq0$ and $\sum_m\pi_m=1$. Even with constant mixing
weights, this prior depends on the observed trajectory. The final estimate is
\begin{equation}
 \widehat c=\underset{c_j\geq c_{\min}}{\arg\min}
 \left\{\tfrac1{18}\|\overline{\mathsf H}c-\bar r\|_2^2+
 \tfrac12(c-c_{\rm pr})^\top\Lambda(c-c_{\rm pr})\right\},
 \label{eq:inverse}
\end{equation}
with $\Lambda_{jj}=s_K\exp(\beta_j)(\eta/\eta_{\rm ref})^2(N_{\rm ref}/N)^2$.
Precision denotes penalty strength. This empirical scaling increases regularization with noise
and reduces it with denser sampling. The $\eta^2/N^2$ factor follows the variance of an
average of $N^2$ independent pixels as a guide; smoothing and dependent features preclude
an exact inverse-covariance interpretation. The reference scales and bounded log-precisions
are specified in \ref{app:training}. They are not an inverse-covariance formula.
For clean observations, $\Lambda=0$. Enumeration of the 32 free-coordinate
sets finds the feasible minimum in double precision, with gradients through the selected candidate
when its active set remains unchanged locally.

Keeping this constrained solve fixed lets three regularization branches test the value of the
final stage. The global branch fits
one mixture and five precisions across training laws. A non-neural smoothness branch tunes
a second-difference penalty on validation laws. SCCT adapts the prior and precision
to each trajectory. Reduced images supply context,
while the physical moments retain the observation resolution. Moment-only and nested maps
isolate the inputs and outputs that contribute. Every branch uses the same admissible family
and constrained solve, without separate neural optimization for each trajectory; architecture and training
details are in \ref{app:training}.

\subsection{Residual control and inverse stability}\label{sec:conditional}
The two stages meet in the stability of the realized inverse: correction
reduces the moment residual, while regularization controls its amplification.
The following bound makes their roles explicit.
\begin{proposition}[Conditional stability of the constrained inverse]\label{prop:stability}
For given observed moments, prior $c_{\rm pr}$, and precision $\Lambda$, let $c_\star$ be
feasible, and set $\rho:=g-Kc_\star$. If $K+\Lambda\succ0$, then Eq.~\eqref{eq:inverse} has a unique
minimizer $\widehat c$ and
\begin{equation}
 \|\widehat c-c_\star\|_2\leq
 \frac{\|\rho\|_2+\|\Lambda(c_{\rm pr}-c_\star)\|_2}
 {\lambda_{\min}(K+\Lambda)}.
 \label{eq:bound}
\end{equation}
\end{proposition}
\noindent The proof tests the minimizer's variational inequality with $c_\star$ and applies
Cauchy--Schwarz; \ref{app:moments} gives the proof. The numerator separates the moment residual $\rho$ from the
regularized prior displacement, while the denominator measures their amplification. Positive
precision makes $K+\Lambda$ positive definite for noisy observations; at zero noise, full column rank of
$\overline{\mathsf H}$ is required. Uniform pure-phase equilibria fail this condition because $b_j(\pm1)=0$ for every $j$. The bound remains valid with active constraints and organizes moment,
interpolation, numerical-model, and family errors once $\|\rho\|_2$ is controlled independently.
To relate this coefficient bound to the reported force-error metric, define the Gram matrix on the evaluation
grid by
\begin{equation}
 Q_{ij}=a_0^2\sum_\ell b_i(u_\ell)b_j(u_\ell),\qquad
 \|G_{\widehat c}-G_{c_\star}\|_2^2
 =(\widehat c-c_\star)^\top Q(\widehat c-c_\star).
 \label{eq:force_gram_main}
\end{equation}
When $G_\star=G_{c_\star}$, multiplying the right side of Eq.~\eqref{eq:bound} by
$\sqrt{\lambda_{\max}(Q)}/\|G_\star\|_2$ bounds $e_G$. The same coefficient error can
produce different force errors depending on its direction. For an out-of-family truth,
Eq.~\eqref{eq:force_gram_main} measures distance to a chosen admissible reference, and an
approximation discrepancy remains. The proof in \ref{app:moments} concerns a unique optimizer for the realized moments.
Identifying all five force directions additionally requires informative observations.
Physical-law uniqueness over unrestricted nonlinearities and gauge classes is a broader
partial differential equation (PDE) inverse
problem \cite{NiLai2024Uniqueness,BrunkEggerHabrich2023,KianLiimatainenLin2024}.
Here $\rho$ collects the effects of observation noise, temporal reconstruction, spatial
discretization, and force-family discrepancy. Equation~\eqref{eq:bound}
is deterministic for the observed system; a stochastic accuracy guarantee
would also require a bound on these contributions to $\rho$. To see how
regularization filters that residual,
consider an interior solution.

The implemented system then gives
\begin{equation}
\widehat c-c_\star=
\underbrace{(K+\Lambda)^{-1}\rho}_{r_{\rm data}}+
\underbrace{(K+\Lambda)^{-1}\Lambda(c_{\rm pr}-c_\star)}_{r_{\rm prior}}.
\label{eq:error_components}
\end{equation}
An active-set solution adds the corresponding boundary term on constrained coordinates. The original
SCCT and global solutions analyzed below are interior, so Eq.~\eqref{eq:error_components} gives an
exact truth-dependent decomposition for those synthetic tests.

Equation~\eqref{eq:error_components} motivates two separate tests: reducing the
residual by correcting the equation, and adapting how the prior influences the solve.
The local prior and precision responses in \ref{app:training} explain
why additional network flexibility need not improve a data-informed direction.

\section{Experimental design and results}\label{sec:design}
The numerical study tests each link in the path from observations to
transferable force. We first remove
measurement bias, choose a balance for the observation source, assess the
learned regularizer, and reuse the recovered force. All observations are
numerically generated, including examples with specified physical parameters. The
primary trajectory contains 61 frames over $[0,0.012]$ on a $128^2$ periodic grid,
with iid Gaussian noise $\eta\in\{0.001,0.01,0.03\}$ and separate clean controls.
This densely sampled trajectory permits both balances to be tested; reduced-frame studies
then expose their sampling dependence. At fixed horizon, 61, 31, and 21 frames correspond
to observed intervals $2$, $4$, and $6\times10^{-4}$, respectively; holding the
generator and law fixed separates observation-time quadrature from simulation error.
The noise levels probe the tradeoff between the leading $O(\eta^2)$ nonlinear-substitution
bias of an unadjusted polynomial and the variance of its Hermite-adjusted counterpart;
they are test conditions, not universal transition points. Positive forces and initial shapes are sampled
independently. The original collection has 768 training laws and separate 160-law validation
and test sets, with degrees one, two, and four represented in every split. Their distributions,
all later independent test sets, and selection/reuse rules appear in
Appendix Table~\ref{app:tab:study_map}.

Neural maps minimize relative squared force error on synthetic training labels. Global
mixtures are trained on the same laws, numerical penalties use validation,
and trained models remain unchanged unless a study explicitly retrains them. An original
checkpoint is a validation-selected model held unchanged in later tests. The final SCCT
checkpoints were selected by force loss; the later degree-supervision comparison retained
the force-only starting checkpoint for every seed (\ref{app:training}).
The trajectory comparator is our periodic, constant-$q$, positive-family adaptation of
\cite{Cohen2026Differentiable}. \ref{app:forecast} reports inference inputs and costs.

Evaluation measures the entire force on a uniform phase grid,
\[
 e_G=\left(\frac{\sum_\ell|G_{\widehat c}(u_\ell)-G_\star(u_\ell)|^2}
 {\sum_\ell|G_\star(u_\ell)|^2}\right)^{1/2}.
\]
Percent errors are $100e_G$, and differences use percentage points (pp).
A constitutive law is the independent unit; realizations, seeds, and geometries are grouped
within law. Original benchmark means average seed-specific errors; stated coefficient
ensembles average coefficients before evaluating the norm. Paired bootstrap intervals
resample laws: 95\% for secondary comparisons and Bonferroni $1-0.05/4=98.75\%$ for each locked
study's four primary contrasts. Failures are retained and validation fixes alarm thresholds.
\ref{app:protocol} specifies aggregation, resampling, and calibration data.

\begin{figure}[!htbp]\centering
\begin{subfigure}[t]{0.25\textwidth}\centering\vspace{0pt}
\includegraphics[width=\linewidth]{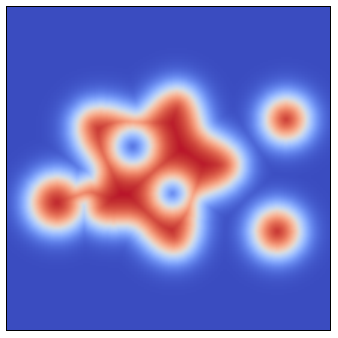}
\caption{Initial perforated field}
\end{subfigure}\hspace{0.5in}
\begin{subfigure}[t]{0.25\textwidth}\centering\vspace{0pt}
\includegraphics[width=\linewidth]{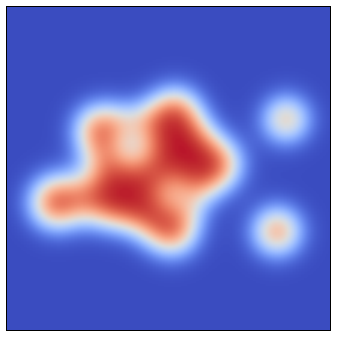}
\caption{Reference, $t=0.004$}
\end{subfigure}\hspace{0.5in}
\begin{subfigure}[t]{0.25\textwidth}\centering\vspace{0pt}
\includegraphics[width=\linewidth]{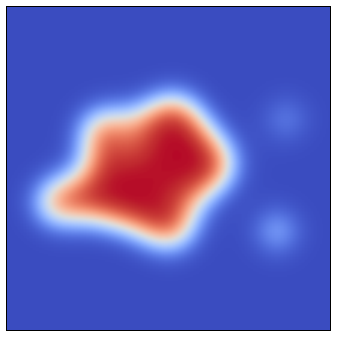}
\caption{Reference, $t=0.008$}
\end{subfigure}\par\smallskip
\begin{subfigure}[t]{0.25\textwidth}\centering\vspace{2pt}
\includegraphics[width=\linewidth]{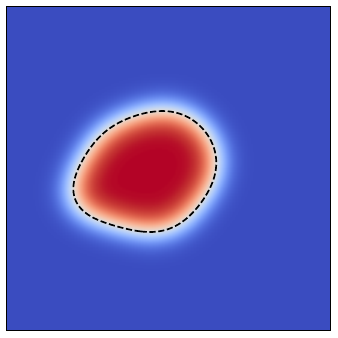}
\caption{SCCT, $t=0.020$}
\end{subfigure}\hspace{0.5in}
\begin{subfigure}[t]{0.25\textwidth}\centering\vspace{2pt}
\includegraphics[width=\linewidth]{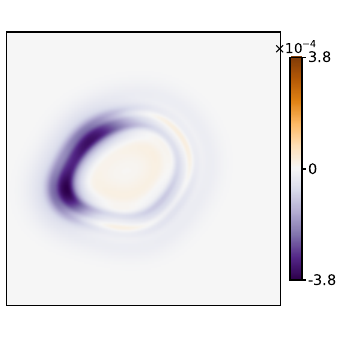}
\caption{SCCT minus reference}
\end{subfigure}\hspace{0.5in}
\begin{subfigure}[t]{0.25\textwidth}\centering\vspace{2pt}
\includegraphics[width=\linewidth]{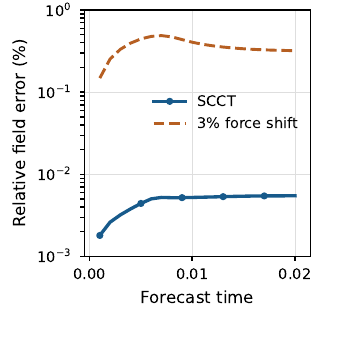}
\caption{Forecast error over time}
\end{subfigure}
\caption{Transfer to an unseen topology-rich geometry. Blue and red in
(a)--(d) denote the negative and positive phase. The dashed black line in
(d) is the exact-law zero contour. All fields use the same exact initial
condition on a $512^2$ grid and Fourier Crank--Nicolson/second-order Adams--Bashforth
(CNAB2) solver \cite{AscherRuuthWetton1995}; only the force differs. The previously inferred
SCCT force comes from the separate noisy wobble trajectory and has
$100e_G=0.096\%$ on the phase grid. Panel (f) compares
$100e_u(t)=100\|u_{\mathrm{model}}(t)-u_\star(t)\|_2/\|u_\star(t)\|_2$
with the $3\%$ force-rescaling control on a logarithmic error axis. The
signed difference in (e) uses units of $10^{-4}$.}
\label{fig:physical_case}
\end{figure}

The results follow the path in Figure~\ref{fig:workflow} from noisy observations
to an inferred force and a field forecast. To make that final step tangible,
Figure~\ref{fig:physical_case} transfers the previously inferred three-seed SCCT coefficient
mean from a noisy wobble trajectory ($\eta=0.01$) to a new perforated,
five-lobed geometry with a narrow neck and two satellite regions. The inverse
is not refitted on this field. The exact-law reference begins with three
positive components and two holes; both holes close by $t=0.004$, and only
one component remains at $t=0.008$. At $t=0.020$, the transferred force gives
$0.0055\%$ relative field error, compared with $0.321\%$ after a deliberate
$3\%$ rescaling of the exact force. Halving the reference time step changes
its final field by $0.00061\%$. Comparing the exact-law $256^2$ and $512^2$
solutions from the same sampled initial field changes the final field by
$0.012\%$. This single-law illustration precedes the
population-level inverse and forecast tests below.

\begin{figure}[!htbp]\centering
\begin{subfigure}[t]{0.25\textwidth}\centering\vspace{0pt}
\includegraphics[width=\linewidth]{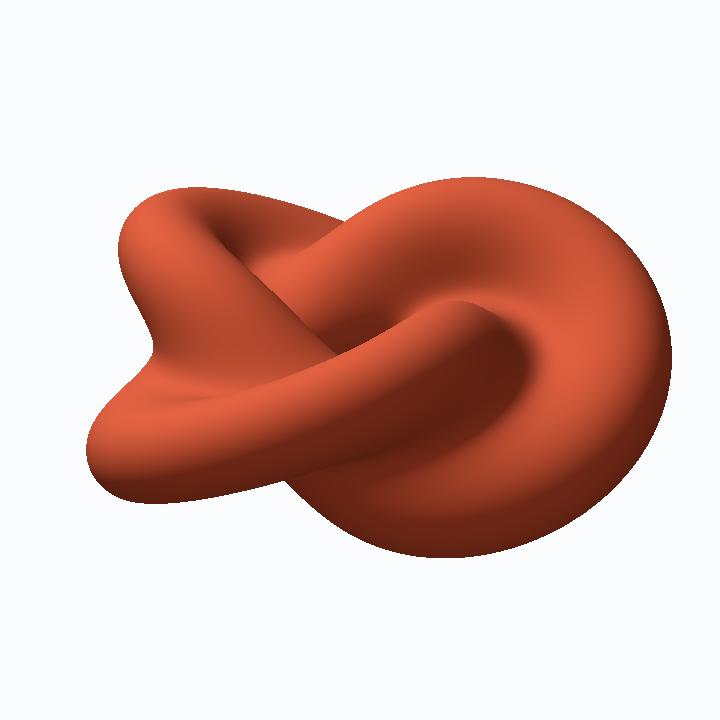}
\caption{Trefoil tube, degree four}
\end{subfigure}\hspace{0.5in}
\begin{subfigure}[t]{0.25\textwidth}\centering\vspace{0pt}
\includegraphics[width=\linewidth]{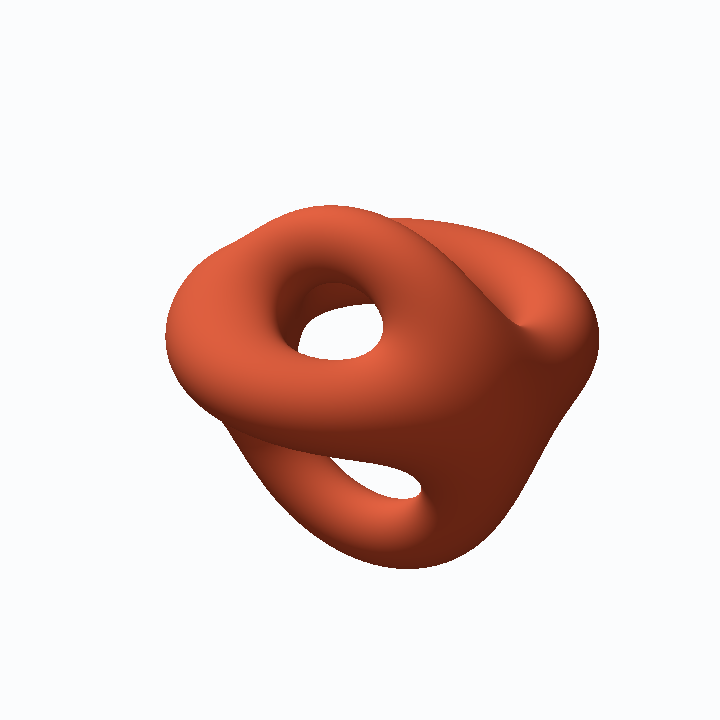}
\caption{Orthogonal loops, degree one}
\end{subfigure}\hspace{0.5in}
\begin{subfigure}[t]{0.25\textwidth}\centering\vspace{0pt}
\includegraphics[width=\linewidth]{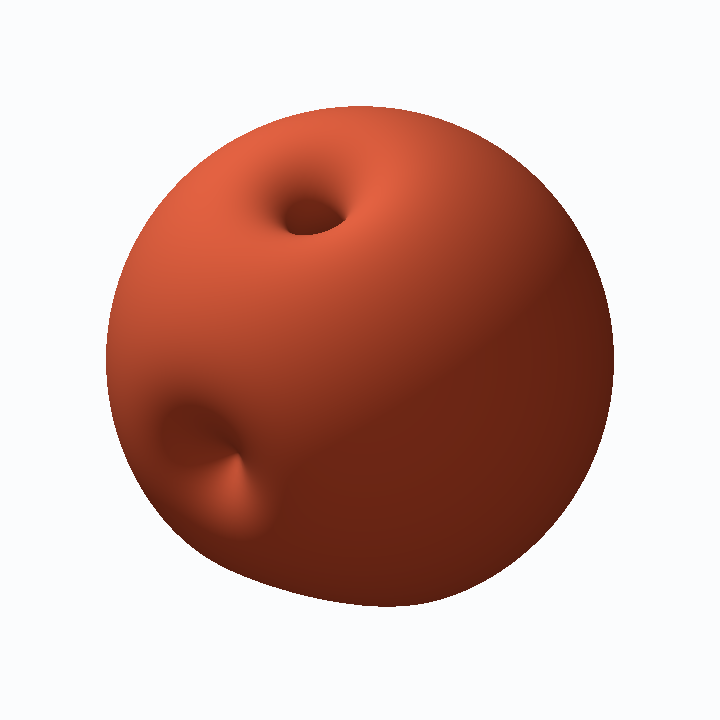}
\caption{Perforated shell, degree two}
\end{subfigure}\par\smallskip
\begin{subfigure}[t]{0.25\textwidth}\centering\vspace{0pt}
\includegraphics[width=\linewidth]{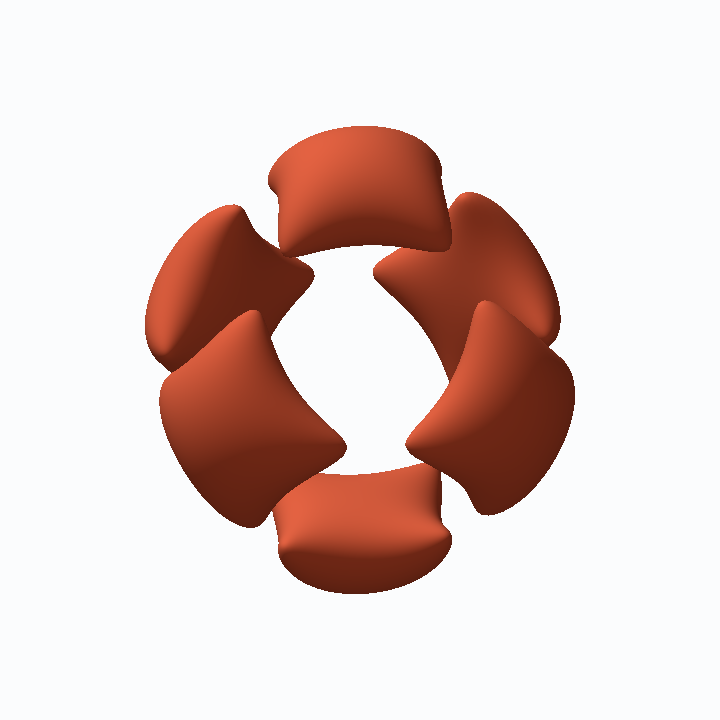}
\caption{Reference after breakup}
\end{subfigure}\hspace{0.5in}
\begin{subfigure}[t]{0.25\textwidth}\centering\vspace{0pt}
\includegraphics[width=\linewidth]{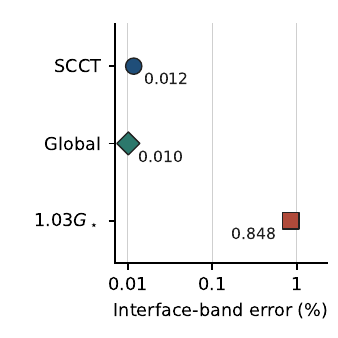}
\caption{Same-time force control}
\end{subfigure}\hspace{0.5in}
\begin{subfigure}[t]{0.25\textwidth}\centering\vspace{0pt}
\includegraphics[width=\linewidth]{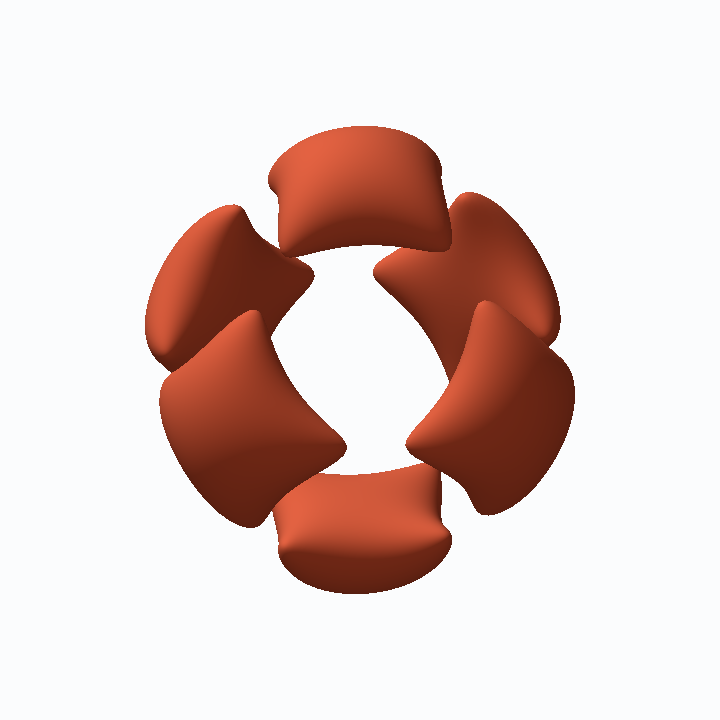}
\caption{SCCT after breakup}
\end{subfigure}
\caption{Three-dimensional transfer and a force-sensitivity control.
Panels (a)--(c) render the SCCT $u=0$ isosurface at $t=0.004$ for one
law of each degree on a different geometry. The same three laws are reused
on the seven-handle ring cage, not counted as additional independent laws:
(d) and (f) compare the exact and SCCT isosurfaces at the same time,
$t=0.008$, for the degree-four law. Both have six positive components.
Panel (e) compares the SCCT and global inverses with a $3\%$ rescaling
of the exact force on that same ring cage at $t=0.008$; its horizontal
error axis is logarithmic. Each comparison starts from the same exact
initial field: the $128^3$ grid is used in (a)--(c), and the
$192^3$ grid in (d)--(f). Interface error is
$100\|\chi_B(u_{\rm model}-u_\star)\|_2/\|\chi_Bu_\star\|_2$,
where $B=\{|u_\star|<0.9\}$. No three-dimensional data enter the inverse.}
\label{fig:three_dimensional}
\end{figure}

The inferred scalar force also enters the three-dimensional equation
$u_t=\Delta_3u-G(u)$ without changing its coefficients. To test that transfer
beyond one illustrative law, we continue twelve laws spanning three force
orders, inferred from independent two-dimensional trajectories, on three unseen
three-dimensional geometries (Figure~\ref{fig:three_dimensional}). At
$t=0.012$, mean diffuse-interface errors are $0.0215\%$ for SCCT and
$0.0229\%$ for the global-mixture inverse. Their paired difference is
$-0.00145$ pp, with a degree-stratified law-bootstrap $95\%$ interval
$[-0.00877,0.00350]$ pp; SCCT has lower error on five of twelve laws.
SCCT has the lower mean for degrees one and two, whereas the global mixture
has the lower mean for degree four. Across all three orders and geometries,
both two-dimensional inverses give accurate three-dimensional fields without
three-dimensional observations or refitting; this experiment does not rank
them reliably.

A ring cage probes a more complicated geometry while separating geometric
motion from force sensitivity. For a balanced profile $u\simeq Q(d)$ with
$Q''=G(Q)$, the leading normal-profile diffusion and reaction terms
cancel in
$u_t=\Delta u-G(u)$, leaving curvature-driven interface motion
\cite{AllenCahn1979}. The exact seven-handle body splits into six components
between saved times $0.007$ and $0.008$. SCCT, the global inverse, and even
a deliberately $3\%$ rescaled exact force share this zero-level topology
sequence. Yet on the same $192^3$ grid at $t=0.008$, the representative
SCCT forecast has $0.0117\%$ diffuse-interface error versus $0.848\%$ for
the rescaled-force control. Across the three reused laws, SCCT and global
mean errors are $0.0158\%$ and $0.0154\%$. Thus field accuracy, not the
breakup label alone, supports transfer of the inferred force.

\subsection{Correcting the equation matters more than retraining}\label{sec:factorial}
The first comparison asks whether retraining can compensate for an
uncorrected equation. We cross Hermite correction with
independent views, rebuild the features, reselect numerical penalties, and refit every
global and adaptive branch. At $\eta=0.03$, the complete estimator is most accurate on
160 laws (Table~\ref{tab:conventional_correction}). Removing views increases adaptive-map
error by $0.24$ pp $[0.18,0.31]$; removing polynomial correction increases it by
$1.40$ pp $[1.24,1.56]$. These paired analyses average seed-specific differences within
law, whereas the table reports coefficient-ensemble means.

\begin{figure}[!htbp]\centering
\includegraphics[width=0.49\textwidth]{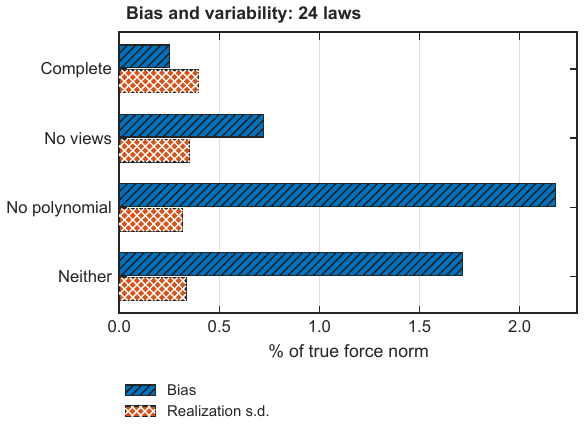}
\caption{
  Bias and realization spread after retraining on a separate 24-law set with five realizations at
  $\eta=0.03$; spread uses divisor five. Diagonal and cross hatching distinguish
  bias and realization spread in grayscale.
}
\label{fig:correction_sampling}
\end{figure}

The high-noise penalty raises a second question: does it reflect displacement or variability?
Across 24 laws with five realizations, removing polynomial correction increases the
finite-repeat force displacement by nearly an order of magnitude without increasing
the spread (Figure~\ref{fig:correction_sampling}). The effect persists after retraining.
The nonmonotone factorial follows from doubled view variance: without polynomial
correction, view generation amplifies substitution bias. \ref{app:corrections}
defines the finite-repeat diagnostics.

This reduction in displacement has a variance cost. At $\eta=0.001$, the
uncorrected smoothness branch gives $0.031\%$ error versus $0.057\%$ complete, while
the retrained complete map reaches $0.027\%$. All constructions coincide on clean data.
Across the tested levels, both corrections become beneficial from $\eta=0.01$ upward.
The benefit of correction depends on the size of the nonlinear noise bias.

The non-neural smoothness inverse tests whether this ranking depends on the learned map.
It gives the same moderate-noise ordering: removing views,
polynomial correction, or both increases its error by $0.440$ pp $[0.330,0.564]$,
$2.727$ pp $[2.351,3.129]$, and $1.613$ pp $[1.357,1.888]$, respectively
(95\% law-bootstrap intervals). Agreement across the two regularizers locates the
improvement in the estimating equation. After removing the noise bias, we next
test whether the temporal balance follows the forward generator.

\begin{table}[!htbp]\centering\footnotesize\setlength{\tabcolsep}{3pt}
\caption{
  Noise-correction factorial on 160 laws. Values are mean force errors (\%). Smoothness is a
  non-neural positive inverse with a second-difference penalty selected on validation; the adaptive
  map is retrained for each construction. All constructions coincide on clean data. The $\times$
  column is each row's $\eta=.03$ mean divided by that of the full noise-adjusted construction.
}\label{tab:conventional_correction}
\begin{adjustbox}{max width=\textwidth}
\begin{tabular}{@{}lS[table-format=1.3]S[table-format=1.3]S[table-format=1.3]cS[table-format=1.3]S[table-format=1.3]S[table-format=1.3]c@{}}\toprule
\rowcolor{gray!12}Construction & \multicolumn{4}{c}{Smoothness inverse} & \multicolumn{4}{c}{Retrained map}\\
\cmidrule(lr){2-5}\cmidrule(l){6-9}
\rowcolor{gray!12} & {$\eta=.001$} & {$\eta=.01$} & {$\eta=.03$} & {$\times$} & {$\eta=.001$} & {$\eta=.01$} & {$\eta=.03$} & {$\times$}\\\midrule
Full noise adjustment & 0.057 & 0.264 & 0.490 & \cellcolor{gray!10}$1.00\times$ & 0.027 & 0.198 & 0.416 & \cellcolor{gray!10}$1.00\times$\\
Without views & 0.123 & 0.345 & 0.930 & \cellcolor{gray!10}$1.90\times$ & 0.030 & 0.291 & 0.656 & \cellcolor{gray!10}$1.58\times$\\
Without polynomial correction & 0.149 & 0.966 & 3.218 & \cellcolor{gray!10}$6.57\times$ & 0.098 & 0.697 & 1.807 & \cellcolor{gray!10}$4.34\times$\\
Neither correction & 0.031 & 0.712 & 2.103 & \cellcolor{gray!10}$4.29\times$ & 0.056 & 0.362 & 1.391 & \cellcolor{gray!10}$3.34\times$\\\bottomrule
\end{tabular}
\end{adjustbox}
\end{table}

\subsection{The observation source selects the estimating balance}\label{sec:unknown}
Noise adjustment alone does not address numerical-model error: the temporal
balance must still describe the
source of the frames. We cross two generators with two balances to test that requirement,
then place the resulting weak moments alongside public implementations.

For a source distinct from the IMEX generator, twenty-four laws are evolved with a refined
fourth-order exponential time-differencing Runge--Kutta (ETDRK4) solver
\cite{CoxMatthews2002,KassamTrefethen2005} and sampled onto the inference grid. Time and space
refinement change the final forward field negligibly at the stated precision. Identical arrays are
supplied to all methods; resolution checks are in \ref{app:corrections}.

\paragraph{Changing the temporal balance}
Table~\ref{tab:observation_source} crosses the two generators with the two estimating balances. At
61 frames and $\eta=0.01$, replacing an assumed IMEX balance with observation-time weak moments
reduces the original SCCT checkpoint's error from $0.371\%$ to $0.219\%$, with a clean floor reduced from $0.326\%$
to $0.043\%$. The paired path change is $-0.1522$ pp $[-0.1937,-0.1135]$ for SCCT and $-0.1591$ pp
$[-0.2041,-0.1151]$ for global regularization.

\begin{figure}[!htbp]\centering
\includegraphics[width=0.43\textwidth]{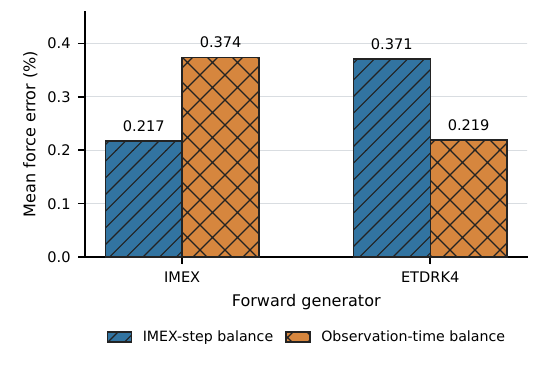}
\caption{Generator-by-balance crossover for SCCT at $\eta=0.01$ on 24 laws.
Each pair of hatched bars uses the same observations; lower mean force error is
better. The preferred balance reverses between IMEX- and ETDRK4-generated trajectories. The
smoothness inverse follows the same ordering (Table~\ref{tab:observation_source}).}
\label{fig:source_crossover}
\end{figure}

\begin{figure}[!htbp]\centering
\begin{subfigure}[t]{0.385\textwidth}\centering\vspace{0pt}
\includegraphics[width=\linewidth]{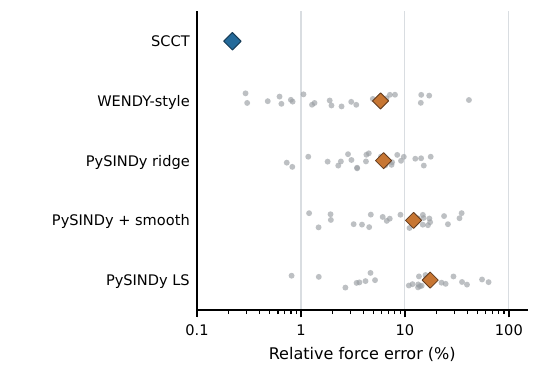}
\caption{$\eta=0.01$}
\end{subfigure}\hspace{0.5in}
\begin{subfigure}[t]{0.385\textwidth}\centering\vspace{0pt}
\includegraphics[width=\linewidth]{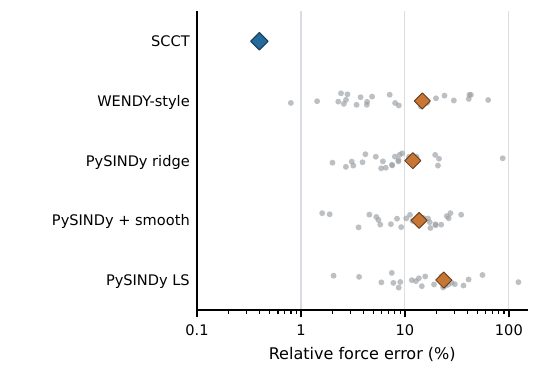}
\caption{$\eta=0.03$}
\end{subfigure}
\caption{Force error on the same 24 ETDRK4-generated laws. Gray points show
the realization-averaged error for each law under the previously evaluated weak baselines;
diamonds show method means on a logarithmic axis. The SCCT diamond is the
corresponding mean in Table~\ref{tab:observation_source}, while baseline
law-level points come from the saved v24 outputs.}
\label{fig:weak_distribution}
\end{figure}

If the source rather than the balance changes, the preference should reverse. For
IMEX-generated observations, observation-time moments give a
$0.336\%$ clean floor and $0.374\%$ SCCT error at $\eta=0.01$, versus $0.217\%$ with a discrete
balance derived from the IMEX steps. Among these two constructions, the discrete balance ranks first for IMEX
observations and the observation-time weak balance ranks first for ETDRK4 observations. In
Table~\ref{tab:observation_source}, the balance chosen for each source has lower error than the alternative in
every column, at each noise level, on clean observations, and for the neural as well as the smoothness
inverse. The reversal persists across these two generators, noise levels, and regularizers
(\ref{app:corrections}). Figure~\ref{fig:source_crossover} shows the crossover, and
Figure~\ref{fig:weak_distribution} shows the spread of external weak-baseline
errors on the same laws; Table~\ref{tab:observation_source} gives the exact means.

\paragraph{External weak-form context}
The crossover selects between our two balances, but does not establish how they compare with
an independently implemented weak construction. We use PySINDy's weak-form library
\cite{deSilva2020PySINDy,Kaptanoglu2022PySINDy}, a Python implementation of SINDy,
with the same arrays, five admissible functions, and known Laplacian coefficient.
Plain weak least squares has
clean error below $0.001\%$, but noisy errors of $17.5/23.6\%$ at $\eta=0.01/0.03$.
Using the same positive smoothness inverse, the public weak moments give $12.1/13.7\%$
versus $0.252/0.493\%$ with noise-adjusted moments (Table~\ref{tab:observation_source}).
This gap compares complete moment constructions; the retrained factorial isolates the
two noise corrections. A covariance-weighted PDE adaptation of weak-form estimation
of nonlinear dynamics (WENDy)
\cite{Bortz2023WENDy} does not resolve an improvement over tuned ridge on these observations
(Appendix Table~\ref{tab:covariance_baseline}). It is distinct from the public
pyWENDy and WENDy--MLE implementations for ordinary differential equations
(ODEs) \cite{Rummel2025WENDyMLE}.

\begin{table}[!htbp]\centering\footnotesize
\caption{
  Observation-source matrix for 61 frames on the $128^2$ grid. Entries are mean force errors (\%) on
  24 laws. Smoothness is this paper's positive constrained inverse with a validation-selected
  second-difference penalty. The last row replaces the noise-adjusted moments by PySINDy's weak library
  \cite{deSilva2020PySINDy,Kaptanoglu2022PySINDy} on the same arrays with the same five admissible
  functions and known Laplacian coefficient. Its first group uses the library's least-squares
  inverse; its second uses this paper's positive smoothness inverse, so the contrast changes
  only the moments. Each $\times$ column divides the $\eta=0.01$ mean by the better balance for the same
  generator. Software details and alternative baseline selections are in
  \ref{app:baselines}.
}\label{tab:observation_source}
\setlength{\tabcolsep}{2.5pt}
\begin{adjustbox}{max width=\textwidth}
\begin{tabular}{@{}llS[table-format=2.3]S[table-format=2.3]S[table-format=2.3]cS[table-format=2.3]S[table-format=2.3]c@{}}\toprule
\rowcolor{gray!12}Generator & Estimating balance & \multicolumn{4}{c}{SCCT, or baseline} & \multicolumn{3}{c}{Smoothness}\\
\cmidrule(lr){3-6}\cmidrule(l){7-9}
\rowcolor{gray!12} & & {clean} & {$.01$} & {$.03$} & {$\times$} & {$.01$} & {$.03$} & {$\times$}\\\midrule
IMEX & IMEX discrete & {$<10^{-5}$} & 0.217 & 0.407 & \cellcolor{gray!10}$1.00\times$ & 0.250 & 0.518 & \cellcolor{gray!10}$1.00\times$\\
IMEX & observation-time weak & 0.336 & 0.374 & 0.508 & \cellcolor{gray!10}$1.72\times$ & 0.407 & 0.600 & \cellcolor{gray!10}$1.63\times$\\
ETDRK4 & assumed IMEX & 0.326 & 0.371 & 0.510 & \cellcolor{gray!10}$1.69\times$ & 0.402 & 0.607 & \cellcolor{gray!10}$1.60\times$\\
ETDRK4 & observation-time weak & 0.043 & 0.219 & 0.398 & \cellcolor{gray!10}$1.00\times$ & 0.252 & 0.493 & \cellcolor{gray!10}$1.00\times$\\\midrule
ETDRK4 & PySINDy weak moments & {$<0.001$} & 17.497 & 23.645 & \cellcolor{gray!10}$79.9\times$ & 12.116 & 13.699 & \cellcolor{gray!10}$48.1\times$\\\bottomrule
\end{tabular}
\end{adjustbox}
\end{table}

\paragraph{Changing regularization with the same moments}
Holding the weak moments fixed isolates the smaller contribution of regularization:
SCCT--global is $+0.0061$ pp $[-0.0212,+0.0327]$. After both maps are
refitted on weak features, adaptive--global is $-0.0304$ pp $[-0.0513,-0.0129]$, with means
$0.195\%$ and $0.225\%$. These exploratory 95\% intervals resample laws after realization averaging
and coefficient ensembling. The transferred model chiefly benefits from the revised temporal
balance; the refitted comparison tests adaptation with that balance held fixed.

The weak path removes internal-stage metadata, but sparse sampling exposes its
temporal-quadrature error in the same-array ETDRK4 study. At 31/21
frames, the weak clean floor rises to $0.470/0.464\%$, and noisy errors at $\eta=0.01$ to
$0.543/0.555\%$. Against the refined trajectory inverse, weak SCCT changes from $-0.115$ pp
$[-0.203,-0.033]$ at 61 frames to $+0.184$ pp $[0.041,0.312]$ at 31 frames. That trajectory control
has a $0.075\%$ clean floor under its finite optimization budget. Pixel averaging and coarse point
sampling also degrade recovery (\ref{app:corrections}); the weak path
removes an internal-stage requirement, but its accuracy still depends on temporal resolution and
the spatial observation operator.

\subsection{Learned regularization is secondary after noise adjustment}\label{sec:decomposition}\label{sec:precision_controls}\label{sec:architecture_locked}
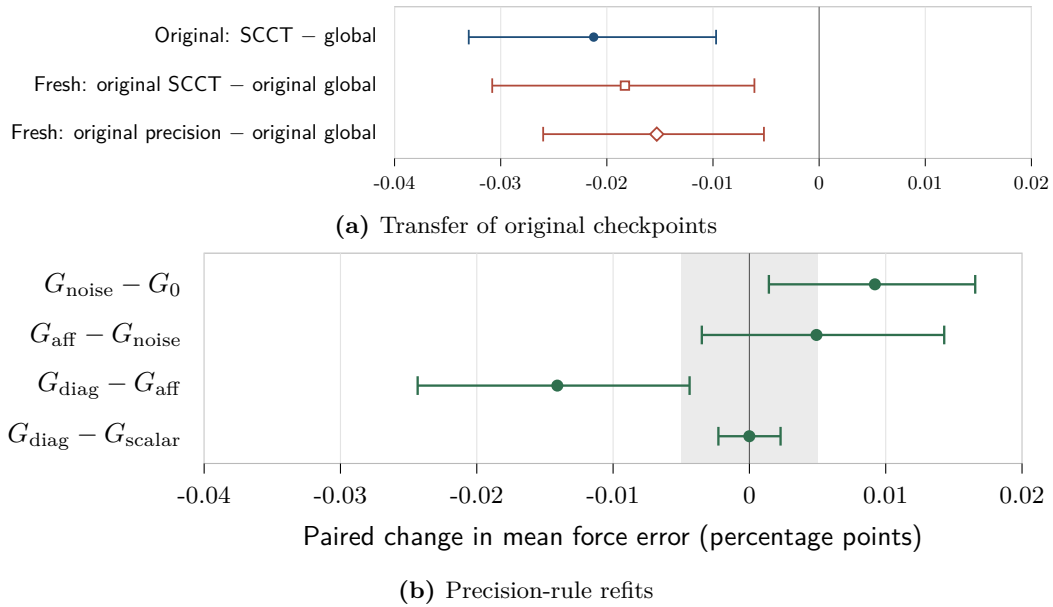
\begin{figure}[!htbp]\centering
\begin{subfigure}[t]{0.79\textwidth}\centering\vspace{0pt}
\resizebox{\linewidth}{!}{%
\definecolor{originalblue}{RGB}{31,78,121}
\definecolor{posthocred}{RGB}{176,74,58}
\begin{tikzpicture}[x=1cm,y=1cm,font=\sffamily\small]
  \def\xmin{0}
  \def\xmax{10.5}
  \def\xzero{7.0}
  \draw[gray!55] (\xmin,3.65) rectangle (\xmax,6.15);
  \foreach \x in {0,1.75,3.5,5.25,7,8.75,10.5} {\draw[gray!22] (\x,3.65)--(\x,6.15);}
  \draw[black!65] (\xzero,3.65)--(\xzero,6.15);
  \node[anchor=east] at (-0.15,5.65) {Original: SCCT $-$ global};
  \node[anchor=east] at (-0.15,4.85) {Fresh: original SCCT $-$ original global};
  \node[anchor=east] at (-0.15,4.05) {Fresh: original precision $-$ original global};
  \draw[originalblue,line width=0.85pt] (1.223,5.65)--(5.301,5.65);
  \draw[originalblue,line width=0.85pt] (1.223,5.53)--(1.223,5.77) (5.301,5.53)--(5.301,5.77);
  \fill[originalblue] (3.283,5.65) circle (2.2pt);
  \draw[posthocred,line width=0.85pt] (1.610,4.85)--(5.933,4.85);
  \draw[posthocred,line width=0.85pt] (1.610,4.73)--(1.610,4.97) (5.933,4.73)--(5.933,4.97);
  \draw[posthocred,line width=0.85pt,fill=white] (3.728,4.78) rectangle (3.868,4.92);
  \draw[posthocred,line width=0.85pt] (2.450,4.05)--(6.090,4.05);
  \draw[posthocred,line width=0.85pt] (2.450,3.93)--(2.450,4.17) (6.090,3.93)--(6.090,4.17);
  \draw[posthocred,fill=white,line width=0.85pt]
    (4.323,4.16)--(4.433,4.05)--(4.323,3.94)--(4.213,4.05)--cycle;
  \foreach \x/\lab in {0/{-0.04},1.75/{-0.03},3.5/{-0.02},5.25/{-0.01},7/{0},8.75/{0.01},10.5/{0.02}} {
    \draw (\x,3.65)--(\x,3.52);
    \node[anchor=north,font=\footnotesize] at (\x,3.48) {\lab};
  }
\end{tikzpicture}
}
\caption{Transfer of original checkpoints}
\end{subfigure}\par\smallskip
\begin{subfigure}[t]{0.79\textwidth}\centering\vspace{0pt}
\resizebox{\linewidth}{!}{%
\definecolor{precisiongreen}{RGB}{47,111,78}
\begin{tikzpicture}[x=1cm,y=1cm,font=\sffamily\small]
  \def\xmin{0}
  \def\xmax{10.5}
  \def\xzero{7.0}
  \fill[gray!16] (6.125,0.25) rectangle (7.875,2.95);
  \draw[gray!55] (\xmin,0.25) rectangle (\xmax,2.95);
  \foreach \x in {0,1.75,3.5,5.25,7,8.75,10.5} {\draw[gray!22] (\x,0.25)--(\x,2.95);}
  \draw[black!65] (\xzero,0.25)--(\xzero,2.95);
  \node[anchor=east] at (-0.15,2.55) {$G_{\mathrm{noise}}-G_0$};
  \node[anchor=east] at (-0.15,1.90) {$G_{\mathrm{aff}}-G_{\mathrm{noise}}$};
  \node[anchor=east] at (-0.15,1.25) {$G_{\mathrm{diag}}-G_{\mathrm{aff}}$};
  \node[anchor=east] at (-0.15,0.60) {$G_{\mathrm{diag}}-G_{\mathrm{scalar}}$};
  \draw[precisiongreen,line width=0.85pt] (7.250,2.55)--(9.899,2.55);
  \draw[precisiongreen,line width=0.85pt] (7.250,2.43)--(7.250,2.67) (9.899,2.43)--(9.899,2.67);
  \fill[precisiongreen] (8.611,2.55) circle (2.2pt);
  \draw[precisiongreen,line width=0.85pt] (6.389,1.90)--(9.501,1.90);
  \draw[precisiongreen,line width=0.85pt] (6.389,1.78)--(6.389,2.02) (9.501,1.78)--(9.501,2.02);
  \fill[precisiongreen] (7.861,1.90) circle (2.2pt);
  \draw[precisiongreen,line width=0.85pt] (2.741,1.25)--(6.232,1.25);
  \draw[precisiongreen,line width=0.85pt] (2.741,1.13)--(2.741,1.37) (6.232,1.13)--(6.232,1.37);
  \fill[precisiongreen] (4.534,1.25) circle (2.2pt);
  \draw[precisiongreen,line width=0.85pt] (6.602,0.60)--(7.400,0.60);
  \draw[precisiongreen,line width=0.85pt] (6.602,0.48)--(6.602,0.72) (7.400,0.48)--(7.400,0.72);
  \fill[precisiongreen] (6.999,0.60) circle (2.2pt);
  \foreach \x/\lab in {0/{-0.04},1.75/{-0.03},3.5/{-0.02},5.25/{-0.01},7/{0},8.75/{0.01},10.5/{0.02}} {
    \draw (\x,0.25)--(\x,0.12);
    \node[anchor=north,font=\footnotesize] at (\x,0.08) {\lab};
  }
  \node[anchor=north] at (5.25,-0.40) {Paired change in mean force error (percentage points)};
\end{tikzpicture}
}
\caption{Precision-rule refits}
\end{subfigure}
\caption{
  Paired mean force-error differences in percentage points; negative values favor the first
  named model. (a) The original contrast uses a 95\% interval; fresh original-checkpoint contrasts
  are post hoc with 98.75\% intervals. (b) Precision-rule contrasts are prespecified with
  Bonferroni 98.75\% intervals. Shading marks the $\pm0.005$ pp effect-size convention,
  not an equivalence margin. Laws are the resampling units; studies are not pooled.
}
\label{fig:paired_effects}
\end{figure}

Only after forming the noise-adjusted equation do we ask how much additional accuracy comes
from observation-dependent regularization. The original benchmark, fresh precision study, and locked structural controls address successively
stricter comparisons (Figure~\ref{fig:paired_effects}; \ref{app:precision}).
The weak-path refit in Section~\ref{sec:unknown} tests adaptation on a different moment construction.

\paragraph{Original benchmark}
On the original 160 test laws, SCCT reduces mean force error from the global mixture's $0.237\%$
to $0.216\%$, a paired difference of $-0.021$ pp with 95\% interval $[-0.033,-0.010]$. Nested
controls designed after this benchmark associate the improvement with precision adaptation rather than expert
mixing: adapting precision alone reproduces the SCCT mean, and adding weight adaptation changes it
by only $-0.003$ pp $[-0.006,0.001]$. The difference is concentrated, with 16 laws selected post hoc from
the same test errors carrying $82\%$ of it, and a truth-dependent decomposition of
Eq.~\eqref{eq:error_components} associates it with a smaller prior-error term (
\ref{app:precision}). Neither analysis is a confirmatory test.

\paragraph{Fresh precision comparison}
Models, selection rules, and the evaluation protocol were fixed before 160
fresh \emph{test} laws (32/32/96 of degrees 1/2/4; three noise levels, two realizations, 960
observations) were generated; the training and validation collections are those of the original
benchmark and are reused here. All controls share one trained global mixing profile and observation-dependent prior
and differ only in the precision map: global ($G_0$), noise-specific ($G_{\rm noise}$), affine in six
observation-only diagnostics ($G_{\rm aff}$), or a 16-unit neural map with one broadcast
($G_{\rm scalar}$) or five ($G_{\rm diag}$) outputs. Mean errors are $0.2465\%$, $0.2557\%$,
$0.2606\%$, $0.2465\%$, and $0.2465\%$: the global, scalar, and diagonal maps agree to four
decimals. The four primary contrasts were prespecified with Bonferroni 98.75\% coverage:
$G_{\rm noise}-G_0=+0.0092$ pp $[+0.0014,+0.0166]$,
$G_{\rm aff}-G_{\rm noise}=+0.0049$ pp $[-0.0035,+0.0143]$,
$G_{\rm diag}-G_{\rm aff}=-0.0141$ pp $[-0.0243,-0.0044]$, and
$G_{\rm diag}-G_{\rm scalar}=-0.000007$ pp $[-0.0023,+0.0023]$.
Noise-specific precision is worse than global precision, while the diagonal map improves on the
affine diagnostic rule. The unresolved contrasts leave small differences undetermined.
Precision-rule refitting is only one test of adaptation. The original checkpoints provide a
separate transfer test on the same 960 observations.
The original and refitted global models both score $0.2465\%$. The original SCCT scores $0.2281\%$ and the
original precision-only map $0.2311\%$, giving post-hoc differences from the original global model of
$-0.0183$ pp $[-0.0308,-0.0061]$ and $-0.0153$ pp $[-0.0260,-0.0052]$ at $98.75\%$ coverage.
Thus the original models retain a small advantage on fresh laws. In contrast, the refitted
diagonal map differs from the original global model by $+0.00003$ pp $[-0.0039,+0.0044]$.

\paragraph{Structural controls}
To test whether greater architectural flexibility explains the transferred advantage,
two sequential studies with distinct locked 96-law test sets compare expert blending, adaptive
gating, and token attention against pooling; all eight Bonferroni 98.75\% intervals include zero
(Appendix Table~\ref{app:tab:structural}). A truth-informed
simplex-grid oracle leaves $0.0522$ pp of mixing headroom with precision held fixed, of which the
learned gate realizes $0.00044$ pp, so candidate complementarity exists, but the tested training protocols exploit little of it. For this force family and tested budgets, the global regularizer is a
strong low-complexity baseline.

Together these results separate transfer of an existing checkpoint from the
effect of a precision rule under common training. The original adaptive model
retains its advantage on fresh laws, but re-fitted precision maps and structural
variants do not establish one under the common protocols. The original and
re-fitted models differ in prior, curriculum, and optimization budget; their
contrast cannot isolate a training mechanism. The prespecified comparisons
above isolate precision rules, whereas the original-checkpoint contrasts are
exploratory.

\subsection{Information loss and model departures limit recovery}\label{sec:detection}
Correction removes identifiable measurement biases, but it cannot create
missing phase information or repair a wrong force family. Phase coverage,
the noise model, and the
constitutive family determine whether the remaining inverse problem is informative.

\begin{figure}[!htbp]\centering
\begin{subfigure}[t]{0.35\textwidth}\centering\vspace{0pt}
\includegraphics[width=\linewidth,trim=0bp 0bp 191bp 17bp,clip]{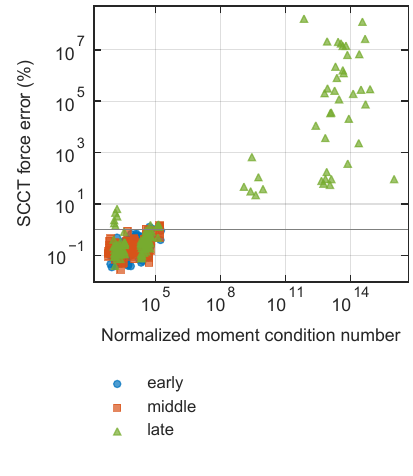}
\caption{Loss of interface information}
\end{subfigure}\hspace{0.5in}
\begin{subfigure}[t]{0.35\textwidth}\centering\vspace{0pt}
\includegraphics[width=0.966\linewidth,trim=198bp 0bp 0bp 17bp,clip]{figures/failure_legible.pdf}
\caption{Forces outside the assumed family}
\end{subfigure}
\caption{
  Different sources of unreliable recovery. (a) SCCT force error versus normalized moment condition
  number for 24 laws, three widths, three windows, and two realizations at $\eta=0.01$; repeated
  trajectories are not independent laws. (b) Mean force error for nonpolynomial forces with 12 laws per
  amplitude, compared with the nonnegative projection benchmark and numerical inverses. The $C=0$
  exponential family is distinct from the separate in-family reference group. The diagnostic
  summary appears in \ref{app:limits}.
}
\label{fig:failure}
\end{figure}

\paragraph{Insufficient phase information}
Twenty-four laws with three initial widths are observed in early, middle, or late 61-frame windows, with
two realizations at $\eta=0.01$. In each width group, 14 of 48 late trajectories have no interface:
moment condition numbers are $1.2\times10^9$--$10^{16}$ and errors at least $22\%$. Trajectories
retaining interfaces have condition numbers below $1.4\times10^5$ and errors below $6.5\%$
(Figure~\ref{fig:failure}a). These ranges describe the tested trajectories; clean-field
occupancy is used for analysis, not inference. For initial width $0.11$, the median error rises from $0.188\%$ in the
early window to $0.395\%$ in the late window, where $29\%$ of trajectories exceed $10\%$ error. Across
all 432 trajectories, the Spearman correlation of log error with log condition number is $0.61$. These
repeated observations come from 24 laws, so the association is descriptive at the law level.

When one trajectory loses phase information, an independent trajectory can broaden the
observed phase range. A separate eight-law experiment with two realizations combines independent initial fields sharing
each law. At $\eta=0.01$, global-inverse mean errors decrease from $0.208\%$ with one $128^2$
trajectory to $0.128\%$ with two and $0.091\%$ with four (\ref{app:limits}).
Independent trajectories offer more phase information when one observation window is weak.
This experiment also increases the observation budget; it does not measure the effect
of geometry diversity at fixed cost.

\paragraph{Sensitivity to the noise model}
Even with adequate phase coverage, correction depends on the measurement-noise model.
On 24 laws with five realizations, temporal median absolute deviation (MAD) estimation maintains high-noise accuracy, while
supplying $0.8\eta$ causes a large error increase (Table~\ref{tab:failure}). On 12 laws with three
realizations, spatial correlation causes much larger errors even when the marginal noise level is
accurately estimated. For measurement covariance $\Sigma$, iid artificial perturbations give
\[
 \operatorname{Cov}(y^+,y^-\mid u)=\Sigma-\eta^2I.
\]
Estimating the marginal variance does not restore independent views: off-diagonal covariance remains
in the estimating equation even if its diagonal is correctly calibrated. Temporal-correlation and
Laplace controls are in \ref{app:limits}.

\begin{table}[!htbp]\centering\small
\caption{
  Noise sensitivity and failure diagnostics. Noise entries are mean SCCT force errors on the
  stated law sets; no ratios are formed across the distinct 24- and 12-law designs.
  Diagnostic thresholds use in-family validation trajectories. AUROC denotes the area under the receiver operating characteristic curve.
}
\label{tab:failure}
\begin{tabular}{@{}p{7.2cm}rr@{}}
\toprule
\rowcolor{gray!12}A: Noise setting, $\eta=0.03$ & Mean $e_G$ (\%) & Laws\\
\midrule
Independent Gaussian, true level / MAD & 0.43 / 0.42 & 24\\
Independent Gaussian, supplied $0.8\eta$ / $1.2\eta$ & 3.72 / 1.08 & 24\\
Spatial correlation, 1 / 2 cells, true marginal s.d. & 6.15 / 8.55 & 12\\
\midrule
\rowcolor{gray!12}B: Diagnostic (84 laws) & AUROC & Missed / 43\\
\midrule
Degree-two/four expert disagreement & 0.998 & 0\\
Moment-equation residual & 0.939 & 6\\
Final-block forecast residual & 0.580 & 38\\
Smallest singular value & 0.500 & 41\\
\multicolumn{3}{@{}p{11.5cm}@{}}{Disagreement and moment residual each flag 1/12 in-family trajectories.
Disagreement also flags 7/12 accurately recovered extrapolation laws.}\\
\bottomrule
\end{tabular}
\end{table}
\paragraph{Departures from the constitutive family and diagnostic limits}
An accurate noise model cannot make a force outside the assumed family admissible.
To test this separate limitation, consider the nonpolynomial family
\[
G_\star(u)=-200u(1-u^2)\exp[A+B(2u^2-1)+C\sin(2\pi u^2)].
\]
We sample 12 laws
at each of five amplitudes $|C|=0,0.15,0.30,0.45,0.60$. Together with 12 coefficient-extrapolation
and 12 in-family laws, this gives 84 IMEX trajectories at $\eta=0.01$, one per law. The evaluation-only
oracle projects the true force onto the nonnegative coefficient family, relaxing the estimator's
strict positive floor. At $|C|=0.30$, this lower approximation benchmark gives $1.12\%$, SCCT
force error is $2.67\%$, and its distance from the best projection is $2.31\%$ of the true force
norm. These norms are not additive components. The inverse fits trajectory-weighted estimating
moments, whereas the oracle projection minimizes force error over the full phase grid using the true
law. Their objectives differ under misspecification; the excess over projection cannot be assigned
solely to neural regularization. All compared inverses show excess error beyond the best projection
(Figure~\ref{fig:failure}b).

Because the true law is unavailable in applications, we next ask whether observed
diagnostics reveal this departure from the assumed family. Alarms define failure as $e_G>1\%$ and use thresholds set at the 95th percentile of in-family
validation trajectories, targeting a 5\% validation false-alarm rate for each score.
Expert disagreement detects all 43 failures, but also flags one of 12 in-family
trajectories and seven of 12 accurately recovered extrapolation laws (Table~\ref{tab:failure}). Small
held-out forecast residuals miss most failures. The smallest-singular-value score is uninformative
in this out-of-family study, with AUROC 0.500 and 41 of 43 failures missed. Conditioning
explains the late-window experiment but fails as a general test of family adequacy. Across this
study, every observation-only score either misses out-of-family laws or also alarms on accurately
recovered laws. The tested diagnostics cannot certify an estimate from an unknown law; developing a
reliable certificate remains an open problem.
\subsection{Recovered laws support reuse within that boundary}\label{sec:new_forecasts}
Within those boundaries, the decisive test is whether an inferred force remains
useful beyond the observed trajectory. We therefore reuse each law from an
unseen initial field and separately ask how an external surface-tension
measurement sets its energy scale.

Twelve new laws, four per degree, are inferred once from 61 IMEX frames at $\eta=0.01$ and used
without refitting from unseen initial fields. Exact-law ETDRK4 references test two planar-front
widths to $t=0.02$ and nine topological cases per law to $t=0.06$, covering annulus closure,
dumbbell pinch-off, and coalescence. Validation selects resolved geometries; exact-law controls are
tilted along one coefficient direction to prescribed force errors.

\begin{figure}[H]\centering
\begin{subfigure}[t]{0.35\textwidth}\centering\vspace{0pt}
\includegraphics[width=\linewidth,trim=0bp 41bp 194bp 16bp,clip]{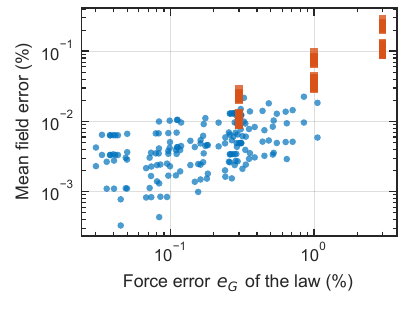}
\caption{Planar-front relaxation}
\end{subfigure}\hspace{0.5in}
\begin{subfigure}[t]{0.35\textwidth}\centering\vspace{0pt}
\includegraphics[width=\linewidth,trim=194bp 41bp 0bp 16bp,clip]{figures/fig_e5_forecast.pdf}
\caption{Topological tests}
\end{subfigure}\par\smallskip
\includegraphics[width=0.78\textwidth,trim=0bp 0bp 0bp 169bp,clip]{figures/fig_e5_forecast.pdf}
\caption{
  Forecast field error versus force error from exact unseen initial fields on 12 laws. Circles pool
  laws inferred by the four methods; squares are exact laws tilted to prescribed force errors along one specified
  direction. (a) Planar-front relaxation; (b) annulus, dumbbell, and two-disc cases. Each law
  generates several dependent forecast cases.
}
\label{fig:new_forecasts}
\end{figure}

On these new initial fields, the inferred laws retain field accuracy. SCCT, moment-only, and global
regularization give similar front and topological field errors, while the finite-budget
trajectory inverse has larger errors (Table~\ref{tab:new_forecasts}). Tilting the exact law
in one coefficient direction increases field error with force error
(Figure~\ref{fig:new_forecasts}). The median field sensitivity per unit force error is
about twice that of SCCT in this imposed direction; \ref{app:forecast} defines
the paired ratio. This is a directional control, not a comparison at equal force-error norm.

These continuous-field errors distinguish the imposed force perturbations more clearly than
event labels do: all four inference methods reproduce 107 of 108 events,
and the tilted laws reproduce all 108. The single disagreement lies at the merging--vanishing
threshold. Force recovery supports field prediction beyond the observed trajectory,
but these saturated event labels and similar field errors under inferred laws provide no
additional architecture-specific advantage.

\begin{table}[H]\centering\small
\caption{
  Forecasts from exact initial fields on 12 laws, with each force inferred once. Field errors are
  time-averaged relative norms (\%); agreement counts nine topological cases per law. The last
  column gives front/event field error relative to SCCT. Event-time resolution is
  $2\times10^{-4}$, so a pooled sub-bin mean does not imply sub-bin resolution for individual
  events. Front-width results are in \ref{app:forecast}; the finite-budget trajectory inverse
  adapts \cite{Cohen2026Differentiable} to the positive force family.
}
\label{tab:new_forecasts}
\begin{adjustbox}{max width=\textwidth}
\begin{tabular}{@{}lrrrrrc@{}}
\toprule
\rowcolor{gray!12}Law used & $100e_G$ (\%) & Front field & Event field & Agreement & Time error & {$\times$SCCT}\\
\midrule
SCCT & 0.201 & 0.0052 & 0.0062 & 107/108 & $1.4\times10^{-4}$ & \cellcolor{gray!10}$1.00/1.00\times$\\
Moment map & 0.187 & 0.0051 & 0.0061 & 107/108 & $1.4\times10^{-4}$ & \cellcolor{gray!10}$0.98/0.98\times$\\
Global mixture & 0.200 & 0.0052 & 0.0061 & 107/108 & $1.4\times10^{-4}$ & \cellcolor{gray!10}$1.00/0.98\times$\\
Trajectory inverse & 0.350 & 0.0082 & 0.0109 & 107/108 & $1.4\times10^{-4}$ & \cellcolor{gray!10}$1.58/1.76\times$\\
Tilted exact, 0.3\% & 0.300 & 0.0170 & 0.0233 & 108/108 & $3.7\times10^{-5}$ & \cellcolor{gray!10}$3.27/3.76\times$\\
Tilted exact, 1\% & 1.000 & 0.0566 & 0.0775 & 108/108 & $1.2\times10^{-4}$ & \cellcolor{gray!10}$10.88/12.50\times$\\
Tilted exact, 3\% & 3.000 & 0.1699 & 0.2315 & 108/108 & $3.1\times10^{-4}$ & \cellcolor{gray!10}$32.67/37.34\times$\\
\bottomrule
\end{tabular}
\end{adjustbox}
\end{table}

\FloatBarrier
\paragraph{External scale and conditional uncertainty}
Forecasts depend on $G$, leaving the energy scale unresolved. With $H'=G$ and $H(\pm1)=0$, an external
positive surface tension $\sigma$ fixes
\begin{equation}
 C_H=\int_{-1}^{1}\sqrt{2H(s)}\dd s,\qquad
 \widehat\eps=\sigma/C_H,\qquad \widehat F=\widehat\eps^{\,2}H.
 \label{eq:calibration}
\end{equation}
A stationary planar profile has energy per unit interface length $\eps C_H$ \cite{Modica1987}; $C_H$
is evaluated by Gauss--Legendre quadrature (\ref{app:moments}). Equation~\eqref{eq:calibration} and the
uniqueness of the calibrated pair it produces are taken from the companion study
\cite{JungKim2026STAC}. Surface tension changes the calibrated energy,
not $G$ or its forecasts. Because the dynamics determine only the combination in $G$, $\widehat\eps$ is a
parameter calibrated by $\sigma$, not a measured interface width.

This separation lets us vary the calibration datum without changing a forecast. With $G$
held constant, a $5\%$ increase in supplied $\sigma$ raises $\widehat\eps$ by $5\%$ and $\widehat F$
by $10.25\%$, with unchanged forecasts. In a separate lognormal experiment, $\sigma_{\rm
obs}=\sigma_\star\exp(sz-s^2/2)$ with $z\sim N(0,1)$ and $s=0.05$ gives mean SCCT potential error
$7.97\%$. This random uncertainty is distinct from the deterministic perturbation identity.

The deterministic scaling does not describe statistical uncertainty in the recovered force.
A 100-repeat parametric bootstrap on 24 in-family laws with two realizations per noise level yields
coverage near the nominal 90\% level for $\eps$ in the designed surface-tension-uncertainty setting,
but lower coverage in other conditions (Appendix Table~\ref{app:tab:coverage}). On eight nonpolynomial laws with
$|C|=0.30$, pointwise force coverage falls to 9\% and no norm-ball interval covers the truth. These
intervals remain conditional on the assumed family and noise model.

\FloatBarrier
\section{Discussion and conclusions}
Taken together, the analysis and experiments show that the path from a noisy
movie to a transferable force begins with the estimating
equation. Hermite features address nonlinear substitution, opposing views
separate instruments from response noise, and a balance based on the observation
source controls clean numerical residuals. The retrained factorial assigns
the largest accuracy change to the moment corrections; the IMEX--ETDRK4
crossover shows why the temporal balance must follow the data generator.
Covariance modeling \cite{Bortz2023WENDy,Rummel2025WENDyMLE} can address
residual dependence beyond these corrections.

After the moments are formed, the conditional bound tracks how their residual
and the prior displacement affect the positive inverse. This explains both
the smaller role of learned regularization in the tested regime and its
dependence on available phase information. The original SCCT checkpoints
retain a small fresh-law advantage, but common-training precision maps and
structural controls perform close to the global regularizer. Their different
training histories preclude assigning that advantage to one architectural
feature. The usefulness of adaptation depends on the inverse and data
distribution \cite{AlbertiDeVito2021OptimalTikhonov,Hertrich2026LearnedRegComparison}.

The recovered force is then tested beyond the observed trajectory. It forecasts
unseen planar fields and three-dimensional evolutions without a new inverse
solve, while external surface tension sets the energy scale independently.
The ring cage makes the distinction between geometry and constitutive accuracy
plain: perturbed forces can share its breakup label while producing much larger
diffuse-interface error. This transfer evidence is conditional on an informative
movie and the assumed force family. Late frames that have lost the interface,
correlated noise, and out-of-family laws challenge those assumptions; the
multi-trajectory and diagnostic tests show how information can be added and
where present alarms remain insufficient. Experimental images, unknown
kinetics, and broader constitutive families are the next settings in which to
test the same sequence of moment construction, inverse analysis, and transfer.

\section*{Declaration of competing interest}
The authors declare that they have no known competing financial interests or personal relationships that could have appeared to influence the work reported in this paper.

\section*{Acknowledgements}
This work was supported by a National Research Foundation of Korea (NRF)
grant funded by the Korean government (MSIT) (No. RS-202625589376).

\appendix
\counterwithin{table}{section}
\renewcommand{\thetable}{\Alph{section}.\arabic{table}}
\section{Study protocol and statistical units}\label{app:protocol}
All comparisons are organized around the generating constitutive law, which
is the independent experimental unit.
Noise realizations, frame counts, initial fields, and training seeds are repeated measurements or
model variations within a law; they are averaged within law before law-level resampling.  Unless a
table states otherwise, the learned predictions shown average coefficients across training seeds
before the force error is evaluated.

\begin{table}[!htbp]\centering\scriptsize
\caption{Study overview. Test counts are independent constitutive laws, not trajectories or frames.}
\label{app:tab:study_map}
\begin{adjustbox}{max width=\textwidth}
\begin{tabular}{@{}>{\raggedright\arraybackslash}p{2.3cm}>{\raggedright\arraybackslash}p{2.6cm}>{\raggedright\arraybackslash}p{4.2cm}>{\raggedright\arraybackslash}p{4.5cm}@{}}
\toprule
Study & Test units & Models and aggregation & Role of selection and reuse\\
\midrule
Original learning comparison & 160 laws; three noise levels & SCCT, neural controls, global and numerical inverses; seed errors averaged within law & Original validation-selected checkpoints; later decomposition is exploratory\\
Locked structural controls & Two disjoint sets of 96 laws; two realizations & Three-seed coefficient ensembles & Choices fixed on 192 training and 48 validation laws before each test set\\
Correction factorial & 160 laws; separate 24-law repeat set & Rebuilt features, refitted global branch, and retrained adaptive maps & Constants and checkpoints reselected on validation for every construction\\
Observation source & 24 laws; three noisy realizations & Discrete and observation-time paths; original and retrained maps & Window support and retrained checkpoints selected on separate validation laws\\
Noise and information & 24 laws; eight-law multi-field subset & Original checkpoints held unchanged and numerical controls on common arrays & Conditions reuse laws and are not independent replications\\
Forecast reuse & 12 new laws; 11 unseen fields per law & One force inferred per law, then reused without another inversion & Forecast geometries chosen using validation laws\\
Out-of-family forces & 84 laws; separate coverage study & Original estimators and parametric bootstrap & Alarm thresholds set by in-family validation trajectories\\
Precision comparison & 160 fresh laws; 960 trajectories & One common mixing profile and five precision rules & Protocol fixed before test generation\\
\bottomrule
\end{tabular}
\end{adjustbox}
\end{table}

The original collection contains 768 training laws, 160 validation laws, and 160 test laws.  Their
degree-one/two/four counts are 128/128/512 in training and 32/32/96 in validation and test.
These proportions define the synthetic target distribution; they are not material-specific priors.  For a
degree-four law,
\[
 c_j=\exp\!\left[A+Bs_j+C(s_j^2-\tfrac12)+z_j\right],\qquad s_j=-1+\tfrac12j,
\]
where $A\sim U[-0.65,0.65]$, $B\sim U[-0.8,0.8]$, $C\sim U[-0.4,0.4]$, and
$z_j\sim N(0,0.07^2)$.  Lower-degree laws use
$c_j^{(m)}=\exp[A+B(-1+2j/m)+z_j]$, $m\in\{1,2\}$ and
$z_j\sim N(0,0.08^2)$, followed by degree elevation.  Initial shapes and laws are sampled
independently.  Shapes include ellipses, angular perturbations, and two-circle unions; profile width
is uniform on $[0.085,0.14]$.  The interfacial parameter is log-uniform on $[0.04,0.1]$.

Synthetic calibration uses $\sigma_\star=\eps_\star C_{H_\star}$.
This generated calibration datum is excluded from the network input.

Force errors use 1001 equally spaced phase values.  Paired intervals resample laws and retain all
within-law grouping.  Multiplicity-adjusted structural and precision studies use Bonferroni
$98.75\%$ intervals for four prespecified contrasts; other intervals are descriptive $95\%$
law-bootstrap intervals. Random streams, split identifiers, and law-level results
allow the numerical comparisons to be repeated.

\section{Moment identities and inverse properties}\label{app:moments}
This appendix gives the algebra behind the admissible force, the moment
corrections, and the conditional inverse bound used in the main text.
\subsection{Admissible family and numerical normalization}
For $0<u<1$, positive Bernstein coefficients give $G_c(u)<0$; odd symmetry gives
$G_c(u)>0$ for $-1<u<0$. Thus $H_c$ increases from the left well to a central maximum and
decreases to the right well, with $H_c\geq0$ and $H_c(\pm1)=0$.
Moreover $G_c'(0)=-a_0c_0$ and $G_c'(\pm1)=2a_0c_4$.
To retain this well structure beyond the phase interval, the implemented extension is
\[
 H_c(u)=a_0c_4(|u|-1)^2,\qquad
 G_c(u)=2a_0c_4\big(u-\operatorname{sgn}u\big),\qquad |u|>1.
\]
It joins the interior primitive with two continuous derivatives and gives a restoring force
outside the wells. Hermite correction uses the polynomial interior basis, without clipping;
overshoots of reconstructed clean stages contribute a model residual rather than a Gaussian
bias identity for this exterior extension.

To implement this family, we use $a_0=200$ and $c_{\min}=10^{-6}$.
The scale normalizes typical coefficients near unity; strict positivity ensures
nondegenerate wells. The floor safeguards the numerical solve.
Surface-tension calibration uses 96-point Gauss--Legendre quadrature. The integrand
$\sqrt{2H_c}$ is generally nonpolynomial, so the quadrature is approximate;
a 96-to-192-node refinement check gave a maximum relative change below
$1.5\times10^{-14}$ in the available numerical check. Quadrature error is separate
from force inference.

Finally, the lower-degree experts must retain admissibility when expressed in the common
five-coefficient family. Degree elevation from $m\in\{1,2,4\}$ to degree four uses
\[
 E^{(m)}_{jk}=\frac{\binom mk\binom{4-m}{j-k}}{\binom4j},
 \qquad \binom ab=0\quad\text{when }b\notin\{0,\ldots,a\}.
\]
Each row sums to one, so elevation preserves both the polynomial and the coefficient floor.
A convex mixture of elevated experts remains admissible.

\subsection{Filters and instrument coverage}
Variance reduction begins with a temporal filter whose normalized weights are
$(1-(m/5)^2)^5$, $m=-5,\ldots,5$; its zero endpoints
leave nine effective taps. The spatial filters are periodic Gaussian convolutions with widths
$N/64$ and $N/32$ grid cells. A preliminary instrument smoother uses $N/128$ cells.
These chosen fractions retain the same physical smoothing supports across grids.
After filtering, the instruments must cover the phases that distinguish coefficients.
In Eq.~\eqref{eq:instruments_main}, $\mu_k=k/8$ and $h_\psi=0.18$ give overlapping coverage
of $u^2\in[0,1]$: adjacent centers are $0.125$ apart, less than the Gaussian width.
The nine instruments overdetermine five coefficients but do not guarantee rank.
The resulting moment rows are put on a common scale: row normalization uses a
$10^{-12}$ denominator floor to avoid division by zero.
For the observation-time path, independent validation selects temporal support $0.002$ for
61/31 frames and $0.003$ for 21 frames; larger support with fewer frames trades localization for
more samples per quadrature window. These are implementation settings, not universal optima.

\subsection{Known-generator and observation-time balances}
The known-generator branch follows the stage placement used to create the observations.
For $q=1$, the IMEX generator makes three implicit--explicit Fourier steps of length
$\delta=\Delta t/3$ per observation interval.  Summing
$(u^{s+1}-u^s)/\delta=L_hu^{s+1}-G_c(u^s)$ gives
\begin{equation}
 \frac{u_{n+1}-u_n}{\Delta t}
 =L_h\!\left(\frac13\sum_{s=1}^{3}u_{n,s}\right)
 -a_0\sum_jc_j\!\left(\frac13\sum_{s=0}^{2}b_j(u_{n,s})\right).
 \label{app:eq:discrete}
\end{equation}
The unobserved internal states are reconstructed from six neighboring frames by degree-five
Lagrange interpolation; only intervals with complete support are used.  Equation
\eqref{app:eq:discrete} is exact for the actual internal states. Substituting interpolated states adds a reconstruction
residual; reproducing the time-stepping scheme does not remove that residual.

To implement that stage-consistent balance on noisy frames, write
$\widetilde y_n(\theta)$ for degree-five
interpolation from six measured frames at interval fraction $\theta$. For one view the features and
diffusion--time response are
\begin{equation}
\begin{aligned}
X_{nj}&=S_wT\left[\frac13\sum_{s=0}^2\mathcal B_j(\widetilde y_n(s/3);v_{n,s})\right],\\
d_n&=a_0^{-1}S_wT\left[L_h\left(\frac13\sum_{s=1}^3\widetilde
y_n(s/3)\right)-\frac{y_{n+1}-y_n}{q\Delta t}\right].
\end{aligned}\label{eq:imex_features_main}
\end{equation}
Here $L_h$ is the Fourier Laplacian. Writing $w_\ell(\theta)$ for the degree-five interpolation
weights at interval fraction $\theta$, the variance carried into the Hermite-adjusted polynomials is
$v_{n,s}=2\eta^2\sum_\ell w_\ell(s/3)^2$ for a recorrupted view, so it depends on $s$ only through
those weights. A finite-support temporal filter $T$ and periodic Gaussian smoother $S_w$ reduce
high-frequency noise before aggregation. Nine instruments cover squared phase values
as defined in Eq.~\eqref{eq:instruments_main}.
The response and force features now share the reconstructed stages. For each instrument view,
preliminary spatial smoothing precedes interpolation.
Average the three preceding internal positions, clip this average to $[-1,1]$, apply $\psi_k$, and
then apply $T$ and $S_w$ to obtain $Z_k$. This clipping concerns instruments only; the noise-adjusted
polynomial arguments remain unclipped. Nine instruments tested against five basis functions produce
the $9\times5$ system. \ref{app:moments} specifies the filters and overlapping phase
supports; their role is variance reduction and phase coverage.

When internal stages are unavailable, the balance instead follows from multiplying the PDE
by a compactly supported temporal test $\phi$ and integrating $u_t$ by parts:
\[
 \int\phi G(u)\,\dd t=\int\phi\Delta u\,\dd t+q^{-1}\int\phi' u\,\dd t.
\]
This gives the observation-time features in Eq.~\eqref{eq:weak_observed};
their instruments average observed windows without internal-stage interpolation.

\subsection{Gaussian polynomial and cross-noise correction}
The first correction removes bias from nonlinear evaluation of noisy phase values.
For $\xi\sim N(0,v)$, use the generating function
\[
 \exp(ty-vt^2/2)=\sum_{k\ge0}Q_k(y;v)t^k/k!.
\]
Its conditional expectation at $y=u+\xi$ is $\exp(tu)$.  Coefficient comparison proves
$\mathbb E Q_k(u+\xi;v)=u^k$ and the recurrence in
Eq.~\eqref{eq:hermite}.  Replacing every monomial of the Bernstein force basis by its $Q_k$
counterpart removes Gaussian nonlinear-substitution bias at the unnormalized-moment
level when the variance matches the true value.  For interpolated views that variance is
$2\eta^2\sum_\ell w_\ell^2$.

The second correction removes noise shared by an instrument and its response. Let the
original measurement error be $e\sim N(0,\eta^2I)$ and introduce independent
$\xi\sim N(0,\eta^2I)$.  The view errors $e^+=e+\xi$ and $e^-=e-\xi$ satisfy
\begin{equation}
 \operatorname{Cov}(e^+,e^-)=0,
 \qquad e^+,e^-\mid u\ \hbox{are independent }N(0,2\eta^2I).
 \label{app:eq:view_independence}
\end{equation}
Cross-pairing an instrument from one view with features and response from the other consequently
removes the instrument--response noise product in expectation.  With correlated error covariance
$\Sigma$ and artificial variance $\tau^2I$, the remaining cross-covariance is
$\Sigma-\tau^2I$; equal marginal variance alone does not restore view independence.  The expectation identities apply before row normalization and inversion, with the true supplied
variance. Plug-in variance, random normalization, constraints, and regularization change the
finite-sample estimator.

Both corrections establish identities for unnormalized moments, not for the completed
inverse. Writing
$e=\bar r-\overline{\mathsf H}c_\star$ for the residual of the normalized system formed below,
\[
 \mathbb E[g-Kc_\star]
 =\tfrac19\big(\mathbb E\overline{\mathsf H}\big)^{\!\top}\mathbb E[e]
 +\tfrac19\textstyle\sum_k
 \operatorname{Cov}\!\big(\overline{\mathsf H}_{k,:}^{\top},\,e_k\big).
\]
The first term represents the effect of random row normalization and the second the coupling between the
normalized design and its residual. Positivity, regularization, and plug-in variance estimation add
further nonlinear steps. The unbiasedness result applies to the noise-adjusted moments, while
Eq.~\eqref{eq:inverse} is analyzed through the conditional bound below.
\ref{app:moments} gives the derivation.

\subsection{Proof of the conditional stability bound}
To prove Proposition~\ref{prop:stability} for the realized constrained solve, write
$A=K+\Lambda$, $d=\widehat c-c_\star$, and
$b=\rho+\Lambda(c_{\rm pr}-c_\star)$.  The variational inequality for the feasible minimizer,
tested with the feasible vector $c_\star$, gives
$\langle A\widehat c-g-\Lambda c_{\rm pr},c_\star-\widehat c\rangle\ge0$, hence
\[
 \langle Ad,d\rangle\le \langle b,d\rangle
 \le \|b\|_2\|d\|_2.
\]
If $A\succ0$, then $\lambda_{\min}(A)\|d\|_2^2\le\|b\|_2\|d\|_2$; the triangle inequality gives
Eq.~\eqref{eq:bound}.  This proof remains valid with active positivity constraints.  It is a
conditional statement for the realized finite moment system, distinct from uniqueness of an
unrestricted physical law.  The force identity in Eq.~\eqref{eq:force_gram_main} then translates
the coefficient bound through the evaluation-grid Gram matrix.

\section{Architecture, training, and local sensitivity}\label{app:training}
To distinguish learned context from the physical inverse, SCCT does not
predict the force directly. It supplies a positive expert mixture and a diagonal
precision to the common constrained inverse.  Physical moments remain at the observation
resolution; only the neural context is averaged to $32^2$.

\begin{table}[!htbp]\centering\small
\caption{Network architecture.}
\begin{tabular}{@{}>{\raggedright\arraybackslash}p{3.1cm}>{\raggedright\arraybackslash}p{10.2cm}@{}}
\toprule
Component & Implemented specification\\
\midrule
Image context & Thirteen frames; five input channels; three circular $3\times3$ convolution layers with 16 channels and dilations $1,2,1$\\
Phase tokens & Nine Gaussian groups on the signed phase interval at each time; 117 width-48 tokens containing pooled features, phase value, mass, and time\\
Constitutive queries & Five queries containing moment correlations, normalized response, auxiliary inverse, noise level, expert coefficients, and residual indicators\\
Attention & Two width-48 transformer encoder layers, four heads, feed-forward width 96, followed by query-to-token cross-attention\\
Physical experts & Positive Bernstein degrees one, two, and four, elevated to five coefficients and mixed by shared softmax weights\\
Output and solve & Five bounded log-precisions and three logits; 32 active sets enumerate the positive quadratic inverse in double precision\\
\bottomrule
\end{tabular}
\end{table}

The token count is $13\times9=117$: each of the 13 context frames contributes
one token per phase group. Five inequality-constrained coefficients have at most
$2^5=32$ active-coordinate patterns, which can be enumerated directly.
Four attention heads divide the width-48 representation into 12 channels per
head. The image and attention widths keep the model small and are held
fixed across the moment-correction comparisons, so changes in the estimating
equation are evaluated with the same network capacity.

The precision references are $\eta_{\rm ref}=0.01$ and $N_{\rm ref}=128$,
with log-precision outputs bounded to $[-24,0]$. These empirical normalizations
and bounds cover weak to moderate penalties while keeping outputs finite.

The signed phase groups have centers evenly spaced from $-0.9$ to $0.9$ and width $0.18$;
they differ from the squared-phase instruments in Eq.~\eqref{eq:instruments_main}.
A zero query channel keeps the original model's input dimension without introducing an untrained
grid-resolution response.

Training uses a 160-epoch degree-four initializer, 80 mixed-degree epochs, and equal-length 40-epoch
force-only and degree-supervised continuations. Adam with decoupled weight decay
(AdamW) \cite{Loshchilov2019AdamW}, cosine learning-rate decay,
gradient clipping, and validation checkpoint selection control optimization rather than physical
model resolution. Starting with degree four exposes the largest admissible force family
before the lower-degree mixture is introduced; force-loss selection follows the
force-error objective. For replay, the schedules run from $(8,4,2)\times10^{-4}$, respectively,
to $2\times10^{-5}$, with batch size 16, weight decay $10^{-5}$, and gradient-norm limit 10.
Seeds are 17, 29, and 43. The force loss is multiplied by $10^4$ for optimization scaling.
The degree-supervised continuation adds $0.1\mathcal L_{\rm degree}$, with inverse-frequency
class weights $(2,2,\tfrac12)$. Checkpoint selection minimizes validation force loss and includes
the starting checkpoint. All three final SCCT seeds retain that force-only checkpoint; the
original checkpoint's improvement cannot be attributed to degree supervision. The final
transformer has 68,980 trainable parameters; mean-aggregation and moment-only controls have 26,164
and 17,476.  Labels, true coefficients, generating degree, and surface tension are unavailable to
the map at inference; the last enters only the downstream calibration.

The local role of expert mixing is also explicit. Hold the observed moments, precision, and active
set unchanged. In the interior,
\begin{equation}
 \delta\widehat c=(K+\Lambda)^{-1}\Lambda\,\delta c_{\rm pr}.
 \label{eq:prior_response_main}
\end{equation}
With active bounds, the same relation restricts the left-hand inverse to free coordinates; the full
$\Lambda\,\delta c_{\rm pr}$ supplies their right-hand side. For a given trajectory, mixing three constant
experts varies the prior in an affine space of dimension at most two. Its effect is further filtered
by the data and precision: when $\Lambda\succ0$, in precision-scaled coordinates the interior
response has eigenvalues
$1/(1+\mu_i)$, where $\mu_i$ are eigenvalues of $\Lambda^{-1/2}K\Lambda^{-1/2}$. Data-informed
directions attenuate prior changes most strongly.

At an interior solution, holding the moments and prior constant while perturbing the precision gives
\begin{equation}
 \delta\widehat c=(K+\Lambda)^{-1}\delta\Lambda\,(c_{\rm pr}-\widehat c).
 \label{eq:precision_response_main}
\end{equation}
Thus precision sensitivity depends on both the prior-to-estimate displacement and the moment system.
Together, Eqs.~\eqref{eq:prior_response_main} and \eqref{eq:precision_response_main} distinguish
moving the prior from changing its directional influence. They are local relations for given
moments; active-set changes use the corresponding piece of the solution map
(\ref{app:training}).

\section{Correction and observation-source controls}\label{app:corrections}
These controls isolate measurement correction from the choice of temporal
balance. The correction factorial rebuilds every feature cache, reselects numerical constants, refits the
global mixture, and retrains the adaptive moment map for each of the four combinations of Hermite
correction and independent views.  On the separate repeat study, with $R=5$ realizations,
\[
 \bar G=R^{-1}\sum_r\widehat G_r,\quad
 b_R=\|\bar G-G_\star\|_2/\|G_\star\|_2,\quad
 s_R^2=\frac{1}{R\|G_\star\|_2^2}\sum_r\|\widehat G_r-\bar G\|_2^2.
\]
For each law, $R^{-1}\sum_r e_{G,r}^2=b_R^2+s_R^2$.  The displayed $b_R$ is a finite-repeat
displacement across the available repeats; the spread uses divisor $R$.

The source experiment uses 24 independent laws generated by ETDRK4 on a $256^2$ grid with internal
step $2.5\times10^{-5}$ and sampled to $128^2$.  Halving the time step changes the final field by at
most $2\times10^{-11}$ in relative norm; increasing the grid to $512^2$ changes it by at most
$2.6\times10^{-7}$.  The same laws are also evolved with IMEX and CNAB2 generators.
All methods receive identical clean or noisy arrays at 61, 31, or 21 observation times.  These
controls establish the generator--balance reversal in Table~\ref{tab:observation_source}; sparse
observations at fewer times additionally expose temporal-quadrature error.

The CNAB2 control uses an independently written force evaluation and its own discrete balance.
At equal resolution, clean errors with the corresponding generator metadata are at the numerical floor, whereas an
assumed IMEX balance leaves a visible clean residual.  Spatial coarsening is checked separately so
that temporal-model error is not conflated with downsampling.  Together these controls isolate time-scheme consistency from spatial observation error.

\section{External weak-form baseline protocol}\label{app:baselines}
The external comparison holds the laws and observations fixed while changing
the weak-form estimator. It uses PySINDy's weak-form library on the same 24 ETDRK4 laws
and noise arrays as the observation-source study.  Its library is
$\{u_{xx}+u_{yy},b_0(u),\ldots,b_4(u)\}$; each $b_j$ is supplied as its exact monomial expansion so
the weak transform acts on the same admissible functions as the proposed estimator. The
Laplacian coefficient is fixed, and only the five constitutive coefficients are solved. No Hermite
correction, independent views, row normalization, or learned prior is added.

Independent validation selected 300 weak subdomains, spatial support $0.25$, temporal support
$0.002$, and polynomial weight degree four. The placement seed is shared across trajectories so
random test-function placement does not confound method comparisons.  Plain, nonnegative, ridge, and positive
second-difference inverses share those weak features.  Ridge strength is selected on independent
validation laws and then held fixed. The WENDy-style row in Table~\ref{tab:covariance_baseline} is a
transparent PDE adaptation of first-order residual-covariance reweighting with 48 disjoint
space--time tests; the public pyWENDy software, which targets ODE parameter estimation, is not
executed. The clean residual of that PDE adaptation exposes its quadrature/model approximation and
appears alongside its noisy errors. Validation and test laws remain disjoint throughout
these comparisons.

\begin{table}[!htbp]\centering\footnotesize
\caption{Seeded external weak-estimator context on the ETDRK4 test laws. Entries are mean force
errors (\%) after averaging realizations within each of 24 laws. All inverse penalties were selected
on the separate validation laws at $\eta=.01$. EIV--IRLS denotes errors-in-variables,
iteratively reweighted least squares; that row is a PDE adaptation of WENDy's
first-order covariance weighting, not the public pyWENDy implementation; its prescribed 48-test
construction has non-negligible clean error.}
\label{tab:covariance_baseline}
\begin{tabular}{@{}lS[table-format=2.3]S[table-format=2.3]S[table-format=2.3]@{}}\toprule
\rowcolor{gray!12}Estimator & {clean} & {$\eta=.01$} & {$\eta=.03$}\\\midrule
PySINDy weak least squares & 0.000 & 17.497 & 23.645\\
PySINDy non-negative least squares & {$<0.001$} & 11.512 & 13.861\\
PySINDy ridge, validation-selected & 1.786 & 6.231 & 11.920\\
WENDy-style PDE EIV--IRLS & 5.554 & 5.830 & 14.704\\
Noise-adjusted weak path, SCCT & 0.043 & 0.219 & 0.398\\\bottomrule
\end{tabular}
\end{table}

\section{Synthetic trajectories with specified physical parameters}\label{app:physical}
To connect the inverse to an energy interpretation, we consider three
trajectories with specified physical parameters. They share one constitutive
law and are demonstrations, not three
independent-law tests.  Each contains 101 fields on a $256^2$ periodic grid; frames 0--60 are used
for inference.  The coordinate transformation
$x_{\rm new}=2x_{\rm old}-1$, $t_{\rm new}=4t_{\rm old}$ maps the domain to $[-1,1)^2$, changes
$\epsilon$ from $0.03$ to $0.06$, and maps $G_{\rm new}=G_{\rm old}/4$ while leaving the potential
and surface tension unchanged.  Their common generating potential and planar profile are
\[
 F_\star(u)=\tfrac14(1-u^2)^2+\tfrac18u^2(1-u^2)^2,\qquad
 Q(z)=\frac{\sqrt2\tanh(\sqrt3z/2)}{\sqrt{3-\tanh^2(\sqrt3z/2)}}.
\]
Applying $Q(d/\epsilon)$ to signed distance creates diffuse initial data for an
evolving two-dimensional solution. The inclusion breakup is a finite-thickness numerical
example whose behavior includes diffuse-interface effects.

The previously trained 17-frame network is applied to five overlapping windows and its positive coefficients
are averaged.  Because it is not retrained on the 61-frame protocol, its larger correction
sensitivity measures transfer of the 17-frame model. The four-way factorial instead
retrains each construction on the 61-frame protocol.
Forecasts for these trajectories use the inferred force in the Fourier solver without learned field correction.
Forecast metrics include relative field error, interface-band root-mean-square error, and positive-phase area
error.  Exact-field restarts isolate constitutive error; noisy observed restarts represent the
operational forecast and also contain initial-state error.

\section{Precision and structural controls}\label{app:precision}
The precision-only comparison holds one global expert mixture fixed and
changes only the precision rule.
Its six observation-only diagnostics are log noise, log condition number of the row-normalized
moment matrix, log raw moment norm, global-inverse residual, degree-two/four expert disagreement,
and an observed interface fraction.  Standardization uses training data only.  Global,
noise-specific, and affine rules use three-start full-batch limited-memory
Broyden--Fletcher--Goldfarb--Shanno (L-BFGS) optimization \cite{Nocedal1980};
scalar and diagonal neural
rules use a 6--16--1 or 6--16--5 tanh map, three training seeds, 150 AdamW epochs, and validation
checkpoint selection.  Every branch uses the same inverse, coefficient floor, precision bounds,
noise/grid scaling, and force loss.

\begin{table}[!htbp]\centering\footnotesize
\caption{Regularization comparisons. Errors are percentages averaged within law. The
original-checkpoint rows on fresh laws have different training histories.}
\begin{adjustbox}{max width=\textwidth}
\begin{tabular}{@{}lrrr@{}}
\toprule
Map & Laws & Mean & 95th percentile\\
\midrule
Original global weights and precision & 160 & 0.237 & 0.517\\
Original adaptive weights, global precision & 160 & 0.227 & 0.461\\
Original global weights, adaptive precision & 160 & 0.216 & 0.458\\
Original adaptive weights and precision & 160 & 0.214 & 0.427\\
Original SCCT & 160 & 0.216 & 0.449\\
\midrule
Fresh global precision & 160 & 0.2465 & 0.4143\\
Fresh noise-specific precision & 160 & 0.2557 & 0.4357\\
Fresh affine diagnostic map & 160 & 0.2606 & 0.4620\\
Fresh scalar / diagonal neural maps & 160 & 0.2465 / 0.2465 & 0.4267 / 0.4173\\
Original SCCT on fresh laws & 160 & 0.2281 & 0.4458\\
\bottomrule
\end{tabular}
\end{adjustbox}
\end{table}

Two subsequent studies use distinct locked 96-law test sets.  The first compares changes to the
existing expert/attention design; the second compares single priors, blending, gating, direct
positive priors, token pooling, and token attention.  Table~\ref{app:tab:structural} preserves all
eight prespecified contrasts in one place.  Every multiplicity-adjusted interval includes zero;
the experiments do not establish that extra expert or attention complexity improves the
stronger baselines.

\begin{table}[!htbp]\centering\scriptsize
\caption{Locked structural contrasts.  Differences are first minus second in percentage points;
intervals have Bonferroni $98.75\%$ coverage within each four-comparison study.}
\label{app:tab:structural}
\begin{adjustbox}{max width=\textwidth}
\begin{tabular}{@{}>{\raggedright\arraybackslash}p{8.0cm}rr@{}}
\toprule
Comparison & Difference & Adjusted interval\\
\midrule
Joint three-expert attention $-$ original SCCT & $-0.008651$ & $[-0.018862,+0.002955]$\\
Joint three-expert attention $-$ longer moment baseline & $+0.006377$ & $[-0.004351,+0.018247]$\\
Joint three-expert attention $-$ corresponding pooling & $-0.008510$ & $[-0.018804,+0.002919]$\\
Joint nine-candidate attention $-$ corresponding pooling & $-0.006369$ & $[-0.028098,+0.004774]$\\
\midrule
Constant blend $-$ validation-selected single degree two & $-0.002890$ & $[-0.007192,+0.001140]$\\
Refitted blend $-$ validation-selected single degree one & $-0.014796$ & $[-0.031639,+0.001545]$\\
Refitted gate $-$ refitted blend & $+0.001032$ & $[-0.008085,+0.009819]$\\
Token attention $-$ token pooling & $-0.008431$ & $[-0.018124,+0.001147]$\\
\bottomrule
\end{tabular}
\end{adjustbox}
\end{table}

With precision held fixed, a truth-informed simplex oracle leaves $0.0522$ pp of mixing headroom, but
the learned gate realizes only $0.00044$ pp.  The oracle uses the unavailable true force and
quantifies candidate complementarity, not learnability.  On the original laws, the truth-dependent
decomposition in Eq.~\eqref{eq:error_components} associates the SCCT improvement with a smaller prior term;
16 post-hoc favorable laws carry most of the mean difference.  Neither result supplies an
observation-only causal explanation.

The planned normal-approximation power at a $-0.005$ pp alternative was about $94\%$ for diagonal
minus scalar and $6$--$13\%$ for the other contrasts; these are design projections, not observed
test power. The effect-size convention controls study sensitivity, not physical significance.

\section{Information loss and departures from the noise and force models}\label{app:limits}
The limiting cases separate loss of phase information from departures in the
measurement and constitutive models. The noise-law study uses the same
observations while changing the supplied variance or the error law.
For independent Gaussian noise, the temporal median-absolute-deviation estimate reproduces the
known-variance result at the main frame counts.  Underestimating the variance is more damaging than
overestimating it because both Hermite correction and precision become too weak.  Spatially
correlated Gaussian noise retains the same marginal variance but violates
Eq.~\eqref{app:eq:view_independence}; temporal correlation and Laplace noise additionally bias the
temporal scale estimate.  These controls delineate assumptions rather than extend the exact
correction beyond Gaussian independent errors.

Information is varied with initial widths $0.05$, $0.11$, and $0.20$ and with early, middle, and
late observation windows. Late trajectories that have lost every interface produce condition numbers
between $1.2\times10^9$ and $10^{16}$ and force errors above $22\%$; trajectories retaining an interface
remain below $1.4\times10^5$ and $6.5\%$.  Combining raw moments from two or four independently
observed initial fields reduces the global error from $0.208\%$ to $0.128\%$ and $0.091\%$ on the
eight-law control.  This combines geometric information and noise replication and is not an
equal-budget proof against higher spatial resolution.

The out-of-family test uses
\[
G_\star(u)=-200u(1-u^2)\exp\!\left[A+B(2u^2-1)+C\sin(2\pi u^2)\right].
\]
The amplitude $|C|$ varies from zero to $0.6$. The oracle minimizes
$\|G_c-G_\star\|_2$ on the evaluation grid with $c_j\geq0$, a relaxation of
$c_j\geq c_{\min}$; this benchmark does not enter inference. Observation-only alarms are calibrated
at the in-family validation $95$th percentile.  Their summary is:

\begin{table}[!htbp]\centering\small
\caption{Observation-only detection of trajectories with force error above $1\%$ on the 84-law family
study.  Missed failures are counted among 43 failures.}
\begin{adjustbox}{max width=\textwidth}
\begin{tabular}{@{}lrrr@{}}
\toprule
Indicator & AUROC & In-family alarms & Missed failures\\
\midrule
Degree-two/four expert disagreement & 0.998 & $8\%$ & 0\\
Moment-equation residual & 0.939 & $8\%$ & 6\\
Training-seed disagreement & 0.846 & $8\%$ & 5\\
Separate-instrument residual & 0.528 & $8\%$ & 40\\
Final-block forecast residual & 0.580 & $0\%$ & 38\\
Smallest singular value & 0.500 & $25\%$ & 41\\
\bottomrule
\end{tabular}
\end{adjustbox}
\end{table}

Expert disagreement detects every tested failure but also alarms on some accurately recovered laws
outside the training coefficient range.  The moment-equation residual is less sensitive.  These
diagnostics warn about specific failure mechanisms; their false alarms and missed failures
measure the limits of model checking.

\section{Forecast, calibration uncertainty, and cost}\label{app:forecast}
Transfer is evaluated without a second constitutive fit. Twelve laws, four
of each generating degree, are inferred once from 61 IMEX-generated frames at
$\eta=0.01$.  Without another inversion, each recovered law is evolved from two perturbed planar fronts and
nine topological initial fields spanning annulus closure, dumbbell pinch-off, and disc coalescence.
Reference trajectories use the exact force and ETDRK4. The event-time metric is evaluated on
saved frames; one near-threshold coalescence contributes most of the identified-law timing error.
Mean relative front-width errors (\%) for SCCT, moment map, global mixture, and trajectory inverse
are $0.017$, $0.017$, $0.016$, and $0.014$, respectively; the $0.3/1/3\%$ tilted-force controls
give $0.059/0.194/0.566$. The width metric does not rank the identified methods in
the same order as force or field error. Noisy restarts add an initial-state contribution:
at $\eta=0.01$, sharp/broad front field errors rise to approximately $0.035/0.081\%$.  Positive controls tilt the exact
coefficient vector to prescribed force errors, confirming that field and front-width metrics
respond to constitutive perturbations even when event type saturates.

The three-dimensional transfer takes four evenly spaced laws from each
degree-one, degree-two, and degree-four group in the independent 160-law
two-dimensional test set, fixed before evaluating any three-dimensional
outcome. For each noisy trajectory at $\eta=0.01$, the known generating coefficients,
SCCT coefficient ensemble, and global expert solution drive separate forecasts; no
three-dimensional observation enters the inverse. Four laws are assigned to
each of three unseen initial geometries, with all degrees represented on
each geometry: a trefoil tube with satellites,
two orthogonal tubular loops, and a spherical shell with two crossed bores.
The fields have the form $u_0=\tanh(d/0.085)$ for the level-set functions
specified in the accompanying code. The transition from $u=-0.9$ to $0.9$
spans about sixteen cells on the $128^3$ periodic grid, resolving the
diffuse interface rather than only its zero surface. All forecasts use
identical Fourier CNAB2 steps through $t=0.012$. Their paired law-level
errors are independent of the planar forecast averages above. Interface-band
error evaluates voxels where $|u_\star|<0.9$, avoiding domination by the
negative bulk phase. The paired interval resamples four laws within each
force-degree group, retaining the balanced design. For one representative
of each geometry, halving the
time step changes the final exact-law field by at most $0.00081\%$.
After projecting the same initial field onto a $96^3$ grid, the coarse
solution differs from the projected $128^3$ reference by at most
$4.7\times10^{-6}\%$. These comparisons keep the constitutive law fixed.

The ring-cage stress test reuses the three degree-representative laws from
the twelve-law transfer. Its initial positive phase is the union of three
orthogonal toroidal tubes, with signed level set
\[
\begin{aligned}
d(x,y,z)&=\max\bigl\{
a-\sqrt{(\sqrt{x^2+y^2}-R_1)^2+z^2},\\
&\hspace{4.4em}a-\sqrt{(\sqrt{y^2+z^2}-R_2)^2+x^2},\\
&\hspace{4.4em}a-\sqrt{(\sqrt{z^2+x^2}-R_3)^2+y^2}
\bigr\},\\
u_0(x,y,z)&=\tanh\!\bigl(d(x,y,z)/0.085\bigr).
\end{aligned}
\]
Here $(R_1,R_2,R_3)=(0.45,0.46,0.47)$ and $a=0.14$ on the same periodic
cube. The slightly different ring radii avoid exact threefold symmetry;
the tube radius yields a connected body with seven handles and no cavities
on $128^3$, $192^3$, and $256^3$ initial grids. The diffuse transition
between $u=-0.9$ and $0.9$ spans about 24 cells on the $192^3$ grid.
All three coefficient sets for each law are continued there by the same
scheme as above. Topology is computed from positive voxels: connected
components give $b_0$, enclosed negative components give $b_2$, and the
cubical Euler characteristic gives the handle count
$b_1=b_0+b_2-\chi$. The representative degree-four exact, SCCT, and
global-force evolutions retain the zero-level change from one component
and seven handles at $t=0.007$ to six components without handles at
$t=0.008$ on the projected $160^3$ grid and under time-step halving on
$192^3$. The half-step trajectories still have the original zero-level
topology at $t=0.0075$. At $t=0.008$, the relative field difference under
time-step halving is at most $0.00077\%$ for the two inferred forces;
the corresponding projected-grid difference is at most
$1.6\times10^{-5}\%$. These are numerical checks of the sampled
zero-level event, not sub-frame event-time estimates.

Topology depends on the phase threshold inside a diffuse interface. For
the representative inferred solutions at $192^3$, the $u=0.1$ level
has already separated at $t=0.007$, whereas the $u=-0.1$ level is still
connected at $t=0.008$. Deliberate rescaling of the exact force by
$3\%$ retains the same sampled $u=0$ event at $192^3$, but raises its
$t=0.008$ interface-band error to $0.848\%$, compared with $0.0117\%$
for SCCT and $0.0101\%$ for the global inverse on this law. In a
$128^3$ sensitivity screen, $10\%$ and $30\%$ rescalings also retain
the sampled event while increasing interface-band error to $2.79\%$
and $8.04\%$. The zero-level topology therefore checks a geometric
trajectory, whereas the diffuse-interface metric distinguishes forces
in this example. For a planar profile $Q''=G(Q)$, the approximation
$\Delta Q(d)-G(Q(d))\simeq \kappa Q'(d)$ explains the leading curvature
motion; finite-width profiles can still depend on the force.

For each exact-start case, directional sensitivity is
\[
 R_{\rm sens}=\frac{(e_{\rm field}/e_G)_{\rm tilted,\,1\%}}
 {(e_{\rm field}/e_G)_{\rm SCCT}},
\]
where $e_{\rm field}$ is the time-averaged relative field error. The median per-case
ratios are $2.1$ for fronts and $2.0$ for topological cases; ratios of pooled means
are $2.19$ and $2.51$. These compare estimation errors with one imposed direction,
not equal-norm force perturbations. All four inference methods reproduce 107 of 108
event labels and all tilted laws reproduce 108; one near-threshold coalescence
distinguishes merging from vanishing. Event agreement is saturated in
this design and supplies little method discrimination.

The calibration datum obeys the exact transformation
$\widehat\epsilon=\sigma/C_{H_{\widehat c}}$ and
$\widehat F=\widehat\epsilon^2H_{\widehat c}$.  Consequently a relative change in $\sigma$ changes
$\widehat\epsilon$ linearly and the energy quadratically, without changing $G$ or any forecast.
Parametric bootstrap repeats the complete estimator after adding synthetic observation noise and,
when specified, resampling $\sigma$. The synthetic noise is added to a clean trajectory generated
under the inferred law from the spatially smoothed first observation, rather than to the original
noisy movie. The $\epsilon$ intervals are basic bootstrap intervals; pointwise $G$ intervals
use bootstrap percentiles. Table~\ref{app:tab:coverage} separates coverage of force,
energy scale, and out-of-family truth.

\begin{table}[!htbp]\centering\small
\caption{Conditional uncertainty; nominal coverage is $90\%$.  In-family rows use 24 laws and 100
bootstrap repeats per trajectory; out-of-family rows use eight laws.}
\label{app:tab:coverage}
\begin{tabular}{@{}>{\raggedright\arraybackslash}p{8.2cm}rr@{}}
\toprule
Interval and setting & Coverage & Laws\\
\midrule
$\epsilon$, uncertain $\sigma$, $\eta=0.01/0.03$ & $0.90/0.90$ & 24\\
$G$ pointwise, in family, $\eta=0.01/0.03$ & $0.91/0.86$ & 24\\
$\epsilon$, exact $\sigma$, $\eta=0.03$ & $0.73$ & 24\\
$G$ pointwise / norm ball, nonpolynomial $|C|=0.30$ & $0.09/0.00$ & 8\\
\bottomrule
\end{tabular}
\end{table}

For the hardware used in these comparisons, one complete $128^2$, 61-frame evaluation takes approximately $0.118$ s for the
three-network SCCT ensemble, $0.088$ s for the moment map, $0.086$ s for the global or smoothness
inverse, and $2.65$ s for the dense trajectory inverse; moment extraction alone takes about
$0.080$ s. Offline training can be amortized against optimization for each trajectory.
The global and smoothness branches remain faster online and avoid neural training. These timings use the implementation evaluated here on the stated hardware.

\begin{table}[!htbp]\centering\footnotesize\renewcommand{\arraystretch}{1.0}
\caption{
  Online cost on one $128^2$, 61-frame trajectory. Moment methods use observed phase,
  known $q$, and noise level; external $\sigma$ enters only energy calibration.
}
\label{tab:methods}
\begin{adjustbox}{max width=\textwidth}\begin{tabular}{@{}>{\raggedright\arraybackslash}p{2.5cm}>{\raggedright\arraybackslash}p{3.0cm}>{\raggedright\arraybackslash}p{1.8cm}>{\raggedright\arraybackslash}p{2.1cm}r@{}}
\toprule
\rowcolor{gray!12}Method & Constitutive inference & Offline labels & Online optimization & Online (s)\\
\midrule
SCCT & Positive inverse; image/moment map & Force labels & Small inverse & 0.093\\
SCCT ensemble & Three SCCT networks & Force labels & Small inverse & 0.118\\
Moment / nested maps & Positive inverse; moment map & Force labels & Small inverse & 0.088\\
Global mixture & Same inverse, 8 global parameters & Force labels & Small inverse & 0.086\\
Smoothness & Positive inverse, second differences & Validation tuning & Small inverse & 0.086\\
Dense trajectory & Positive family, forward sensitivities & None & Trajectory optimization & 2.65\\
\bottomrule
\end{tabular}\end{adjustbox}
\end{table}

\section*{Data and code availability}
The synthetic data and numerical code supporting the principal results
will be made available by the corresponding author upon reasonable request.

\bibliographystyle{unsrtnat}
\setlength{\bibsep}{0pt}
{\small\bibliography{references}}

\end{document}